 \documentclass[12pt]{amsart}
\usepackage{amssymb}
\usepackage{epsf}
\usepackage[all]{xy}
\numberwithin{equation}{section}

\theoremstyle{plain}
 \newtheorem{theo}{Theorem}[section]
 \newtheorem{prop}[theo]{Proposition}
 
 \newtheorem{coro}[theo]{Corollary}

\theoremstyle{definition}
 \newtheorem{defi}[theo]{Definition}
  
 \newtheorem{exa}[theo]{Example}   
   
 \newtheorem{rem}[theo]{Remark}

\newcommand{\onto}{\twoheadrightarrow} 
\newcommand{\map}[1]{\xrightarrow{#1}}

\newcommand{\inc}{\hookrightarrow}
\newcommand{\abs}[1]{\lvert#1\rvert}    % for absolute value
 \newcommand{\ch}{\mathrm{char}}          % for characteristic
 \newcommand{\id}{\mathit{id}}          % for identity

\newcommand{\comp}{\vDash}              %for compositions (of n)
\newcommand{\parti}{\vdash}              %for partitions (of n)
\newcommand{\rank}{\textrm{rk}}           
\newcommand{\sumsub}[1]{\sum_{\substack{#1}}} % sums with multline subscripts
\newcommand{\Min}{\textsf{Min}}

\newcommand{\cop}{\mathit{cop}}
\newcommand{\odd}{\mathrm{odd}}
\newcommand{\im}{\mathrm{Im}}
\newcommand{\End}{\mathrm{End}}

 \newcommand{\frakS}{\mathfrak{S}}
 
 \newcommand{\calG}{\mathcal{G}}
 
\newcommand{\calN}{\mathcal{N}}
\newcommand{\calP}{\mathcal{P}}
 \newcommand{\calQ}{\mathcal{Q}}
 \newcommand{\calY}{\mathcal{Y}}

\newcommand{\calR}{\mathcal R}
\newcommand{\calE}{\mathcal E}

\newcommand{\OQ}{{\mathcal O}\calQ}
\newcommand{\IQ}{{\mathcal I}\calQ}
\newcommand{\OS}{{\mathcal O}S}
\newcommand{\IS}{{\mathcal I}S}
\newcommand{\ON}{{\mathcal O}\calN}
\newcommand{\IN}{{\mathcal I}\calN}
\newcommand{\OSS}{{\mathcal O}\frakS}
\newcommand{\ISS}{{\mathcal I}\frakS}
\newcommand{\OY}{{\mathcal O}\calY}
\newcommand{\IY}{{\mathcal I}\calY}

\newcommand{\zetaQ}{\zeta_\calQ}
\newcommand{\nuQ}{\nu_\calQ}
\newcommand{\chiQ}{\chi_\calQ}
\newcommand{\epsilonQ}{\epsilon_\calQ}
\newcommand{\BzetaQ}{\Bar{\zeta}_\calQ}
\newcommand{\omegaQ}{\omega_\calQ}

\newcommand{\zetaS}{\zeta_S}
\newcommand{\nuS}{\nu_S}
\newcommand{\chiS}{\chi_S}

\newcommand{\BzetaS}{\Bar{\zeta}_S}

\newcommand{\IP}{{\mathcal I}(P)}
\newcommand{\SSym}{{\mathfrak S}\Sym}   
\newcommand{\YSym}{{\mathcal Y}\Sym}
\newcommand{\QSym}{{\mathcal Q}\Sym}
\newcommand{\NSym}{{\mathcal N}\!\Sym}
\newcommand{\Sym}{{\mathit{Sym}}}
 
 \def\H{{\mathcal H}}

 \def\F{\Bbbk}
 \def\ZZ{{\mathbb Z}}
 \def\QQ{{\mathbb Q}}

\newcommand{\Ch}{{\mathbb X}}
\newcommand{\ECh}{\Ch_+}
\newcommand{\OCh}{\Ch_-}
\newcommand{\ChH}{\Ch(\H)}
\newcommand{\EChH}{\ECh(\H)}
\newcommand{\OChH}{\OCh(\H)}

\begin{document}

 \title[Combinatorial Hopf algebras]{Combinatorial Hopf algebras \\ and
generalized Dehn-Sommerville relations}

 \author[M.~Aguiar, N.~Bergeron, and F.~Sottile]{Marcelo Aguiar, Nantel Bergeron,
 and Frank Sottile}

 \address[Marcelo Aguiar]
 {Department of Mathematics\\ Texas A\&M University\\
 College Station, TX 77843, USA}
 \email{maguiar@math.tamu.edu}
 \urladdr{http://www.math.tamu.edu/$\sim$maguiar}

 \address[Nantel Bergeron]
 {Department of Mathematics and Statistics\\York University\\
 Toronto, Ontario M3J 1P3\\
 CANADA}
 \email{bergeron@mathstat.yorku.ca}
 \urladdr[Nantel Bergeron]{http://www.math.yorku.ca/bergeron}

 \address[Frank Sottile]
 {Department of Mathematics and Statistics\\
 University of Massachusetts\\
 Amherst, MA 01003, USA}
 \email{sottile@math.umass.edu}
 \urladdr{http://www.math.umass.edu/$\sim$sottile}

%16W30 Coalgebras, bialgebras, Hopf algebras
%16W50 Graded rings and modules
%05E05 Symmetric functions
%05A15 Exact enumeration problems, generating functions
%06A11 Algebraic aspects of posets
 \date{\today}
\thanks{Aguiar supported in part by NSF grant DMS-0302423}
 \thanks{Bergeron supported in part by CRC, NSERC and PREA}
 \thanks{Sottile supported in part by NSF CAREER grant DMS-0134860}
 

 %%%%%%%%%%%%%%%%%%%%%%%%%%%%%%%%%%%%%%%%%%%%%%%%%%%%%%%%%%%%%%%%%%
 \begin{abstract}
 A {\em combinatorial Hopf algebra} is a graded connected Hopf algebra
 over a field $\F$ equipped with a character (multiplicative linear
 functional) $\zeta:\H\to \F$.
 We show that the terminal object in the category of combinatorial Hopf
 algebras is the algebra $\QSym$ of quasi-symmetric functions; this explains the ubiquity 
of quasi-symmetric functions as generating  functions in combinatorics. We illustrate this
with several examples.
We prove that every character decomposes uniquely as a product of an {\em even} character and an
{\em odd} character. Correspondingly, every combinatorial Hopf algebra $(\H,\zeta)$ possesses
two canonical Hopf subalgebras on which the character $\zeta$ is even (respectively, odd).
The odd subalgebra is defined by certain 
canonical relations which we call the generalized Dehn-Sommerville relations.  We show that,
for $\H=\QSym$, the generalized Dehn-Sommerville relations are the Bayer-Billera
relations and the 
odd subalgebra  is the peak Hopf algebra of Stembridge. We prove that $\QSym$ is
the product (in the categorical sense) of its even and odd Hopf subalgebras. We also calculate the odd subalgebras of
various  related combinatorial Hopf algebras:  the Malvenuto-Reutenauer Hopf algebra of
 permutations,  the Loday-Ronco Hopf algebra of planar binary trees, the Hopf algebras of symmetric functions
and of non-commutative symmetric functions.
 \end{abstract}

 \maketitle

%%%%%%%%%%%%%%%%%%%%%%%%%%%%%%%%%%%%%%%%%%%%%%%%%%%%%%%%%%%%%%%%%%%%%%%%%%%%%%%%%%
\vspace*{-0.3cm}\begin{figure}[!h]
\parbox{350pt}{\tiny\tableofcontents}
\end{figure}
%%%%%%%%%%%%%%%%%%%%%%%%%%%%%%%%%%%%%%%%%%%%%%%%%%%%%%%%%%%%%%%%%%%%%%%%%%%%%%%%%%

 %%%%%%%%%%%%%%%%%%%%%%%%%%%%%%%%%%%%%%%%%%%%%%%%%%%%%%%%%%%%%%

 \section*{Introduction}

 A {\em combinatorial Hopf algebra} is a pair $(\H,\zeta)$ where
 $\H=\bigoplus_{n\ge 0} \H_n$ is a graded connected Hopf algebra over a
 field $\F$  such that $\dim(\H_n)$ is finite
 for all $n\ge 0$, and $\zeta\colon\H\to \F$ is a character (multiplicative linear functional),
 called its {\em zeta function}.
  A morphism $\alpha: (\H',\zeta')\to (\H,\zeta)$ of combinatorial Hopf algebras
 is a  morphism of graded Hopf algebras such that $\zeta'=\zeta\circ\alpha$.
 The terminal object in the category of combinatorial Hopf algebras
 is the algebra $\QSym$ of quasi-symmetric functions, equipped with a canonical character
$\zetaQ\colon\QSym\to\F$  (Theorem~\ref{T:univ-qsym}). This theory provides a natural framework for
combinatorial invariants encoded via quasi-symmetric generating functions. 

 The algebra $\QSym$ of quasi-symmetric functions was introduced by
Gessel~\cite{Ges} as a source of generating functions for poset partitions (Stanley's $P$-partitions~\cite{Sta72}).
The algebra $\Sym$ of symmetric functions is a subalgebra of $\QSym$.
The graded dual of $\QSym$ is the Hopf algebra $\NSym$ of non-commutative symmetric
functions~\cite{MR, GKal}.

Joni and Rota~\cite{JR} made the fundamental observation that many
discrete structures give rise to natural Hopf algebras whose
comultiplications encode the disassembly of those structures.  This was further developed
by Schmitt~\cite{Sch87,Sch94}.
 A first link between these Hopf algebras and quasi-symmetric functions was found by
 Ehrenborg~\cite{Eh96}, who encoded the  flag vector of a graded poset
as  a  morphism from a Hopf algebra of graded posets to $\QSym$.
 A similar construction was given in~\cite{BS99b}, where it was shown that a
 quasi-symmetric function associated to an edge-labeled
 poset also gives a morphism of Hopf algebras.
 Such quasi-symmetric functions
 encode the structure of the cohomology of a flag manifold as a module
 over the ring of symmetric functions~\cite{BS98,BS_skew}.
 These results were later unified via the notion of Pieri operations on posets~\cite{BMSW,BMSW2}.
 In Examples~\ref{Ex:flagvector}--\ref{Ex:Pieri} we show how all these constructions can be
obtained from Theorem~\ref{T:univ-qsym}  in a very natural manner.

A closely related result to Theorem~\ref{T:univ-qsym} was obtained in~\cite{A}. There, 
the terminal object in the category of {\em infinitesimal} Hopf algebras (equipped with a
multiplicative functional)  was described.
Here, we adapt and expand the constructions of~\cite{A}, including the notions of
eulerian subalgebra (here called {\em odd} subalgebra) and  generalized Dehn-Sommerville
relations.

We review the contents of the paper.

We start by discussing the group of characters of a graded connected Hopf algebra in
Section~\ref{S:characters}. Our first main result (Theorem~\ref{T:decomposition}) states that
any character decomposes uniquely as a product of an even character and an odd character.
Even and odd characters are defined in terms of the involution
$\Bar{h}:=(-1)^{\abs{h}}h$ of a graded connected Hopf algebra (Definition~\ref{D:even-odd-char}).

In Section~\ref{S:combhopf} we introduce the notion of combinatorial Hopf algebras $(\H,\zeta)$.
Associated to the character $\zeta$  are certain canonical characters $\chi$ and $\nu$;
these are closely related to the decomposition of $\zeta$ as a
product of an even and an odd character. We describe some combinatorial Hopf algebras involving
partially ordered sets in Examples~\ref{Ex:Rota}--\ref{Ex:dist-lattice}, and we uncover the
combinatorial meaning of the canonical characters in each case.

 We recall the definitions of the Hopf algebras of quasi-symmetric functions and non-commutative
symmetric functions in Section~\ref{S:Q-N}. This short section sets the notation to be
used later in the paper, but it contains no new results and may be skipped by the reader 
familiar with these notions.

 In Section~\ref{S:terminal} we introduce the character $\zetaQ\colon\QSym\to\F$ and
 obtain  the following universal property, which is our second main result (Theorem~\ref{T:univ-qsym}):
  \begin{quote}
   {\em For any combinatorial Hopf algebra $(\H,\zeta)$, there is a unique morphism
 of combinatorial Hopf algebras $\Psi\colon(\H,\zeta)\to(\QSym,\zetaQ)$.}
  \end{quote}
 Thus, $(\QSym,\zetaQ)$ is the terminal object in the category of combinatorial Hopf algebras.
We also show that the Hopf algebra of symmetric functions, equipped with the restriction
$\zetaS$ of $\zetaQ$, is the terminal object in the category of {\em
cocommutative} combinatorial Hopf algebras (Theorem~\ref{T:univ-sym}). We illustrate these
results with several examples of a combinatorial nature, including the flag vector of posets and
the chromatic symmetric function of graphs (Examples~\ref{Ex:flagvector}--\ref{Ex:Pieri}
and~\ref{Ex:powerzeta}--\ref{Ex:theta}). We also show that the morphism
$\QSym\to\QSym$ corresponding to the character $\nuQ$ is the map introduced by Stembridge
in~\cite{Ste97}. Other examples will be presented in~\cite{A2}.

There are two Hopf subalgebras of $\H$ canonically associated to $\zeta$. They are the largest
subcoalgebras on which $\zeta$ is even or odd. The definitions and basic properties of these
objects are presented in Section~\ref{S:even-odd}. The odd subalgebra is defined by certain
linear relations which we call the generalized Dehn-Sommerville relations for $(\H,\zeta)$.
For a homogeneous element $h\in\H$ these relations are
\[ \bigl(\id\otimes(\Bar{\zeta}-\zeta^{-1})\otimes \id\bigr)\circ\Delta^{(2)}(h)\ =\ 0\,,\ \ \ \]
or equivalently
\[\bigl(\id\otimes(\chi-\epsilon)\otimes \id\bigr)\circ\Delta^{(2)}(h)\ =\ 0\,.\]
 We show that the generalized Dehn-Sommerville
relations for $(\QSym,\zetaQ)$ are precisely the relations of Bayer and Billera~\cite{BaBi} 
(Example~\ref{Ex:BB}).

 The construction of even and odd subalgebras is natural: 
 a morphism $\alpha:(\H',\zeta')\to(\H,\zeta)$ of combinatorial Hopf
 algebras sends the even (odd) subalgebra of $(\H',\zeta')$ to the even (odd)
subalgebra of $(\H,\zeta)$ (Propositions~\ref{P:S-I-functoriality} and~\ref{P:even-odd-basic}).
Following~\cite{A}, we show how this simple fact implies 
the important result of Bayer and Billera that the flag
vector of an eulerian poset satisfies the generalized Dehn-Sommerville relations (of
$(\QSym,\zetaQ)$).
 First, it is well-known that
one may construct a certain combinatorial Hopf algebra $(\calR,\zeta)$ from graded posets so that
$\zeta$ and $\zeta^{-1}$ are the usual zeta  and M\"obius functions of
posets~\cite{JR,Sch94}. We add the observation that all eulerian posets belong to
the odd subalgebra of $(\calR,\zeta)$. By naturality, the flag vector of an eulerian poset
must belong to the odd subalgebra of $(\QSym,\zetaQ)$, and thus satisfy
the generalized Dehn-Sommerville relations.

 In Section~\ref{S:even-odd-QSym} we describe the even and odd subalgebras of $(\QSym,\zetaQ)$ in
explicit terms. These algebras have basis elements indexed by even and odd compositions,
respectively. In particular, we find that the odd subalgebra is precisely
 the peak Hopf algebra of Stembridge~\cite{Ste97}. We show that these subalgebras are the
terminal objects in the categories of even and odd combinatorial Hopf algebras
(Corollary~\ref{C:E-Pi-univ}). We show that the morphisms
$\QSym\to\QSym$ corresponding to the even and odd parts of $\zetaQ$ 
are projections onto
these subalgebras (Proposition~\ref{P:canonicalproj}) and derive our third main
result: that $\QSym$ is the categorical product of its even and odd Hopf 
subalgebras (Theorem~\ref{T:product}). These results confirm the importance
of Stembridge's Hopf algebra and at the same time 
unveil the mystery behind its construction. 

 In Section~\ref{S:even-odd-Sym} we describe the even and odd subalgebras of $(\Sym,\zetaS)$.
These algebras have basis elements indexed by even and odd partitions. 
We find that the odd subalgebra is precisely
the Hopf algebra of Schur $Q$-functions~\cite[III.8.1]{Mac}.

In Section~\ref{S:other} we describe the odd subalgebras of closely related Hopf algebras:
the Hopf algebra of non-commutative symmetric functions,
the Hopf algebra of permutations of Malvenuto and Reutenauer, and the Hopf algebra of planar
binary trees of Loday and Ronco. 

{\bf Assumption.} Most of our results are valid over any
commutative ring $\F$ in which $2$ is invertible, but for simplicity
we work over a base field $\F$ of characteristic different from $2$.
The results in Section~\ref{S:terminal}, however, are valid over any commutative
ring. The assumption on the characteristic is recalled
in the statement of some theorems, when emphasizing it seems necessary.
 %%%%%%%%%%%%%%%%%%%%%%%%%%%%%%%%%%%%%%%%%%%%%%%%%%%%%%%%%%%%%%
 \section{The group of characters of a graded connected Hopf algebra}\label{S:characters}

\begin{defi}\label{D:character}
Let $\H$ be a Hopf algebra over a field $\F$. A {\em character} of $\H$ is a morphism
of algebras $\zeta:\H\to\F$. Thus $\zeta(ab)=\zeta(a)\zeta(b)$ and $\zeta(1)=1$. We also refer to
characters as {\em multiplicative linear functionals}.
\end{defi} 

Recall that the {\em convolution product} of two linear functionals $\varphi,\psi:\H\to\F$ is
\begin{equation}\label{E:convolution}
\H\map{\Delta_{\H}}\H\otimes
\H\map{\varphi\otimes\psi}\F\otimes\F\map{m_{\F}}\F\,,
\end{equation}
where $\Delta_H$ is the coproduct of $\H$ and $m_\F$ is the 
product of the base field. We denote the convolution
product simply by $\varphi\psi$. This turns the full linear dual of $\H$ 
into an algebra. Let $m_\H$ be the product of $\H$.
Sweedler's dual $\H^\circ$ is the largest subspace of the full linear dual such that
$m_{\H}^*(\H^\circ)\subseteq \H^\circ\otimes \H^\circ$. Equipped with the convolution product and
the comultiplication $m_{\H}^*$, $\H^\circ$ is a Hopf algebra. A character $\varphi$ of
$\H$ is precisely a group-like element of $\H^\circ$: $m_{\H}^*(\varphi)=\varphi\otimes\varphi$.

The set $\ChH$ of characters of an arbitrary Hopf algebra
$\H$ is a group under the  convolution product.
 The unit element is the counit
$\epsilon_\H$, and the inverse of a character $\varphi$ is $\varphi\circ S_\H$, where $S_\H$ is
the antipode of $\H$. In fact, using the antipode axiom and the fact that 
$\varphi$ is a morphism of algebras one finds 
  \begin{multline*}
    \varphi(\varphi \circ S_\H)=
 m_{\F}\circ \bigl(\varphi\otimes( \varphi\circ S_\H)\bigr)\circ\Delta_\H\ =\   m_{\F}\circ
(\varphi\otimes \varphi)\circ (\id\otimes S_\H)\circ\Delta_\H\\
    =  \varphi \circ  m_{\H} \circ  (\id\otimes S_\H)\circ \Delta_\H
    \ =\ \varphi\circ   u_\H \circ \epsilon_\H\ =\ u_\F\circ\epsilon_\H\ 
    =\ \epsilon_\H\,,
  \end{multline*}
where $u_\H$ and $u_\F$ are the unit maps of $\H$ and $\F$. Similarly,
$(\varphi \circ S_\H)\varphi=\epsilon_\H$. 

If $\H$ is cocommutative then $\ChH$ is abelian.

A morphism of Hopf algebras
$\alpha:\H'\to\H$ induces a morphism $\alpha^*:\Ch(\H)\to\Ch(\H')$  defined on
a character $\varphi$ of $\H$ by $\alpha^*(\varphi)=\varphi\circ\alpha$. Since $\alpha$ is a
morphism of algebras, $\alpha^*(\varphi)$ is a character of $\H'$, and since $\alpha$ is a
morphism of coalgebras, $\alpha^*$ is a morphism of groups.

Any graded Hopf algebra $\H$ carries a canonical
automorphism  
\[h\mapsto \Bar{h}:=(-1)^n h\]
for homogeneous elements $h\in \H_n$. This is an involution: $\Bar{\Bar{h}}=h$.
Therefore, it induces an involution $\varphi\mapsto\Bar{\varphi}$ on the character group of $\H$,
with 
\[\Bar{\varphi}(h)=(-1)^n\varphi(h) \text{ \ for \ }h\in\H_n\,.\]

\begin{defi}\label{D:even-odd-char}
A character $\varphi$ of a graded Hopf algebra $\H$ is said to be {\em even} if
\[\Bar{\varphi}=\varphi\]
and it is said to be {\em odd} if
\[\Bar{\varphi}=\varphi^{-1}\,.\]
\end{defi}
 
The even characters form a subgroup of $\ChH$ that we denote by $\EChH$. The set $\OChH$ of
odd characters is a subgroup of $\ChH$ if $\H$ is cocommutative, but not in general.

Any morphism of graded Hopf algebras $\alpha:\H'\to\H$ preserves the canonical
involution; hence, the morphism of groups
$\alpha^*:\ChH\to\Ch(\H')$  preserves even and odd characters. 

\begin{exa}\label{Ex:nu-chi} Let $\zeta$ be an arbitrary character of a graded Hopf
algebra $\H$. There is a canonical way to construct an odd character from $\zeta$. Indeed, let
\[\nu:=\Bar{\zeta}^{-1}\zeta\,.\]
As all these operations are defined on $\ChH$, $\nu$ is a character. In addition,
\[\Bar{\nu}=\Bar{\Bar{\zeta}}^{-1}\Bar{\zeta}=\zeta^{-1}\Bar{\zeta}=\nu^{-1}\,,\]
so $\nu$ is odd. This simple construction is of central importance for this work.

Define also
\[\chi:=\Bar{\zeta}\zeta\,.\]
As before, $\chi$ is character, and 
\[\Bar{\chi}=\Bar{\Bar{\zeta}}\Bar{\zeta}=\zeta\Bar{\zeta}\,.\]
If $\H$ is cocommutative then $\zeta\Bar{\zeta}=\Bar{\zeta}\zeta$, so $\chi$ is even. 

In Remark~\ref{R:squares} we interpret these constructions in terms of the group structure of
$\ChH$.
\end{exa}

{}From now on we assume that $\H$ is a graded connected Hopf algebra. Thus, 
$\H=\oplus_{n\geq 0}\H_n$, $\H_0=\F\cdot 1$, and all structure maps preserve the grading.
We also assume that each $\H_n$ is finite dimensional. We reserve the notation $\H^*$ for
the {\em graded dual} of $\H$. Thus, $\H^*=\oplus_{n\geq 0}(\H_n)^*$. It is a Hopf subalgebra
of Sweedler's dual $\H^\circ$. Given a linear functional $\varphi:\H\to\F$ we set
\[\varphi_n:=\varphi|_{\H_n}\in\H^*_n\,.\]
Since the comultiplication of $\H$ preserves the grading, we have
\begin{equation}\label{E:convolution-components}
(\varphi\psi)_n=\sum_{i=0}^n\varphi_i\psi_{n-i}\,.
\end{equation}
Any character $\zeta$ satisfies $\zeta_0=\epsilon$. Note also that, if
$\ch\F\neq 2$, $\zeta$ is even if and only if $\zeta_{2n+1}=0$ for all 
$n\geq 0$. An odd character, on the other hand, may have no zero
components.

We set out to show that any character decomposes uniquely as a product of an even
character and an odd character. 

\begin{prop}\label{P:squareroot}
Let $\H$ be a graded connected Hopf algebra over a field $\F$ of characteristic different from
$2$. Let $\varphi,\psi:\H\to\F$ be two linear functionals such that $\varphi_0=\psi_0=\epsilon$.
\begin{itemize}
\item[(a)] There is a unique linear functional $\rho:\H\to\F$ such that
\begin{equation}\label{E:squareroot}
\varphi=\rho\psi\rho \text{ \ and \ }\rho_0=\epsilon\,.
\end{equation}
\item[(b)] If $\varphi$ and $\psi$ are characters then so is $\rho$.
\item[(c)] If $\Bar{\varphi}=\psi^{-1}$ then $\rho$ is even and if $\Bar{\varphi}=\psi$
then $\rho$ is odd.
\end{itemize}
\end{prop}
\begin{proof} According to~\eqref{E:convolution-components} to construct $\rho$ we must solve the
following equations in $\H^*$, for each $n\geq 1$:
\[\varphi_n\ =\ 2\rho_n+\psi_n+\sumsub{i+j+k=n\\0\leq i,j,k< n}\rho_i\psi_j\rho_k\,.\]
Since $2$ is invertible, these equations have a unique solution $\rho_n$, constructed recursively
 from $\rho_0=\epsilon$. This proves (a).

{}From~\eqref{E:squareroot} we deduce
\[\varphi\otimes\varphi\ =\ (\rho\otimes\rho)(\psi\otimes\psi)(\rho\otimes\rho)\,,\]
an equation between linear functionals on $H\otimes H$. Suppose now that $\varphi$ and $\psi$ are
characters. As explained above, $m_{\H}^*(\varphi)=\varphi\otimes\varphi$ and
$m_{\H}^*(\psi)=\psi\otimes\psi$. Therefore, applying the morphism of algebras $m_{\H}^*$
to~\eqref{E:squareroot} we obtain
\[\varphi\otimes\varphi\ =\ m_{\H}^*(\rho)(\psi\otimes\psi)m_{\H}^*(\rho)\,.\]
{}From the uniqueness of such decompositions established in part (a) (applied to the Hopf
algebra $H\otimes H$) we deduce
\[m_{\H}^*(\rho)\ =\ \rho\otimes\rho\,.\]
Thus, $\rho$ is a character, proving (b).

Finally, suppose that $\Bar{\varphi}=\psi^{-1}$. We have
\[\varphi=\rho\psi\rho\Rightarrow \psi^{-1}=\Bar{\varphi}=\Bar{\rho}\Bar{\psi}\Bar{\rho}
\Rightarrow \psi=(\Bar{\rho})^{-1}(\Bar{\psi})^{-1}(\Bar{\rho})^{-1}
=(\Bar{\rho})^{-1}\varphi(\Bar{\rho})^{-1}\Rightarrow \Bar{\rho}\psi\Bar{\rho}=\varphi\,.\]
By uniqueness, $\rho=\Bar{\rho}$, i.e., $\rho$ is even. The remaining half of (c) can be
derived similarly.
\end{proof}

The following is the main result of this section.

\begin{theo}\label{T:decomposition}
Let $\H$ be a graded connected Hopf algebra over a field $\F$ of characteristic different from
$2$. Every character $\zeta:\H\to\F$ decomposes uniquely as a product of characters
\[\zeta=\zeta_+\zeta_-\]
with $\zeta_+$ even and $\zeta_-$ odd. In particular, the only character that is both even
and odd is the trivial one (the counit).
\end{theo}
\begin{proof} Proposition~\ref{P:squareroot} allows us to define an even character $\zeta_+$
such that
\begin{equation}\label{E:zetaplus}
\Bar{\zeta}=\zeta_+\zeta^{-1}\zeta_+\,.
\end{equation}
Let $\zeta_-:=(\zeta_+)^{-1}\zeta$. By construction, $\zeta=\zeta_+\zeta_-$, and
\[\Bar{\zeta}_-=(\Bar{\zeta}_+)^{-1}\Bar{\zeta}=(\zeta_+)^{-1}(\zeta_+\zeta^{-1}\zeta_+)=
\zeta^{-1}\zeta_+=(\zeta_-)^{-1}\,,\]
so $\zeta_-$ is odd. This shows existence. 

Consider another decomposition $\zeta=\varphi\psi$ with $\varphi$ even and $\psi$ odd. Then,
\[\Bar{\zeta}=\Bar{\varphi}\Bar{\psi}=\varphi\psi^{-1}\,,\]
while
\[\varphi\zeta^{-1}\varphi=\varphi\psi^{-1}\varphi^{-1}\varphi=\varphi\psi^{-1}\,.\]
Thus, $\Bar{\zeta}=\varphi\zeta^{-1}\varphi$. By uniqueness of such 
decompositions (Proposition~\ref{P:squareroot}) we must
have $\varphi=\zeta_+$, and then also $\psi=\varphi^{-1}\zeta=\zeta_-$.
\end{proof}

\begin{rem}\label{R:squares} The canonical characters $\nu$ and $\chi$ of
Example~\ref{Ex:nu-chi} are related to the decomposition of Theorem~\ref{T:decomposition} as
follows:
\[\nu=(\zeta_-)^2 \text{\ and, if $H$ is cocommutative, \ }\chi=(\zeta_+)^2\,.\]
Indeed, taking inverses on both sides of~\eqref{E:zetaplus} we obtain
$\Bar{\zeta}^{-1}=(\zeta_+)^{-1}\zeta(\zeta_+)^{-1}$ and then
\[(\zeta_-)^2=(\zeta_+)^{-1}\zeta(\zeta_+)^{-1}\zeta=\Bar{\zeta}^{-1}\zeta=\nu\,,\]
by definition of $\nu$. If $H$ is cocommutative then~\eqref{E:zetaplus} becomes 
$\Bar{\zeta}=(\zeta_+)^2\zeta^{-1}$, which by definition of $\chi$ says that $\chi=(\zeta_+)^2$.
\end{rem}

Theorem~\ref{T:decomposition} is complemented by the following result.

\begin{prop}\label{P:decomposition} Let $\H$ and $\F$ be as above.
Then $\OChH$  is a set of representatives for both the left and the
right cosets of the subgroup $\EChH$ of $\ChH$. Moreover,
\begin{itemize}
\item[(a)] $\OChH$ is closed under conjugation by elements of $\EChH$.
\item[(b)] If $\varphi,\psi\in\OChH$ and $\varphi\psi=\psi\varphi$, then $\varphi\psi\in\OChH$.
\item[(c)] $\OChH$ is closed under integer powers, in particular under inversion.
\end{itemize}
\end{prop}
\begin{proof} Let us first derive (a)--(c). If $\zeta=\varphi\psi\varphi^{-1}$ with
$\varphi$ even and $\psi$ odd, then
\[\Bar{\zeta}=\Bar{\varphi}\Bar{\psi}(\Bar{\varphi})^{-1}=\varphi\psi^{-1}\varphi^{-1}=\zeta^{-1}\,,\]
so $\zeta$ is odd. This is (a). 

Suppose $\varphi,\psi$ are odd and  $\zeta=\varphi\psi$. Then
\[\Bar{\zeta}=\Bar{\varphi}\Bar{\psi}=\varphi^{-1}\psi^{-1}\,.\]
If $\varphi\psi=\psi\varphi$ then $\Bar{\zeta}=\zeta^{-1}$, proving (b).
Part (c) follows from (b).

Theorem~\ref{T:decomposition} says that $\OChH$ is a set of representatives for the
left cosets of the subgroup $\EChH$. Either (a) or (c) imply that the same holds true for
the right cosets. For instance, $\zeta=\bigl(\zeta_+\zeta_-(\zeta_+)^{-1}\bigr)\zeta_+$
displays $\zeta$ as a product of an odd character on the left with an even one on the right.
\end{proof}

\begin{coro}\label{C:decomposition-cocom}Let $\H$ and $\F$ be as above.
Suppose in addition that $\H$ is cocommutative. Then the group of characters
splits as a direct product of abelian groups
\[\ChH\cong\EChH\times\OChH\,.\]
\end{coro}

\section{Combinatorial Hopf algebras}\label{S:combhopf} 

\begin{defi}\label{D:CHA}
 A {\em combinatorial Hopf algebra} is a pair $(\H,\zeta)$ where
 $\H$ is a graded connected Hopf algebra over a field $\F$,
 each homogeneous component $\H_n$ is finite-dimensional, 
 and $\zeta:\H\to\F$ is a character. 

 A morphism of combinatorial Hopf algebras
 $\alpha:(\H',\zeta')\to(\H,\zeta)$ is a morphism $\alpha\colon\H'\to\H$ of graded Hopf
 algebras such that the following diagram commutes:
 $$
    \raise -20pt\hbox{
       \begin{picture}(140,58)
       \put(13,40){$\H'$}    \put(95,40){$\H$}
       \put(57,0){$\F$}
       \put(57,50){$\alpha$}    \put(24,16){$\zeta'$}
       \put(88,16){$\zeta$}
       \put(30,43){\vector(1,0){60}}
       \put(25,35){\vector(1,-1){27}}
       \put(95,35){\vector(-1,-1){27}}
      \end{picture}} %.
 $$
\end{defi}
  Combinatorial Hopf algebras over a field $\F$, together with their
 morphisms, form the category of combinatorial Hopf algebras over $\F$.
Combinatorial algebras and coalgebras and their morphisms are defined similarly (in these cases,
$\zeta$ is only required to be a linear functional).

We refer to $\zeta$ as the {\em zeta character} of the combinatorial Hopf algebra.
We introduce the following terminology for the canonical characters associated to $\zeta$
(Example~\ref{Ex:nu-chi}):
\begin{align*}
\zeta^{-1} & =\zeta\circ S_\H \text{ is the  {\em M\"obius character},}\\
\chi & =\Bar{\zeta}\zeta \text{ is the {\em Euler character},}\\
\nu & =\Bar{\zeta}^{-1}\zeta \text{ is the {\em odd character}.}
\end{align*}

\begin{exa} \label{Ex:Rota} (Rota's Hopf algebra). Let $\calR$ be the
$\F$-vector space with basis the set of all isomorphism classes of finite graded posets.  We say
that a finite poset
$P$ is graded if it possesses a maximum element $1_P$ and a minimum element $0_P$, and all
maximal chains in $P$ have the same length. The rank of $P$, denoted $\rank(P)$, is the length of
a maximal chain from $0_P$ to $1_P$. $\calR$
is a graded connected Hopf algebra where the degree
of $P$ is its rank, the multiplication is cartesian product of posets
\[P\cdot Q:=P\times Q\,,\]
the unit element is the poset with one element,
the comultiplication is
\[\Delta(P):=\sum_{0_P\leq z\leq 1_P}[0_P,z]\otimes [z,1_P]\,,\]
and the counit is
\[\epsilon(P):=\begin{cases}1 & \text{ if
$0_P=1_P$,}\\ 0 & \text{ if not.} \end{cases}\] 

Here, and in everything that follows, if $x$ and $y$ are two elements of 
a poset $P$ then $[x,y]$ denotes the subposet $\{z\in P\ \mid\ x\leq z\leq y\}$. 

This Hopf algebra originated in the work of Joni and Rota~\cite{JR}.
Variations of this construction have appeared repeatedly in the literature~\cite{AF,Eh96,
Sch94}. Here, we add the perspective of combinatorial Hopf algebras.

Let $\zeta:\calR\to\F$ be $\zeta(P)=1$ for every poset $P$. This is the {\em zeta
function of posets} in the sense of Rota. Clearly,
$\zeta$ is a morphism of algebras, so $(\calR,\zeta)$ is a combinatorial
Hopf algebra.  It follows that the M\"obius character
$\zeta^{-1}$ is the classical M\"obius function $\mu$ of posets, in the sense that 
\[\zeta^{-1}(P)=\mu([0_P,1_P])\]
 is the value of the M\"obius function of the poset $P$ on the interval $[0_P,1_P]$. In fact,
from $\zeta^{-1}\zeta=\epsilon$ we deduce the defining  recursion for the M\"obius function: 
\[\zeta^{-1}([x,x])=1 \text{ for every }x\in P \text{ and }\sum_{x\le z\le y}
\zeta^{-1}([x,z])=0  \text{ for every $x<y$ in $P$.}\]
 This approach to M\"obius functions of posets is due to
Schmitt~\cite{Sch94}.

Consider the Euler character $\chi$. We have
\[\chi(P)=\sum_{z\in P}(-1)^{\rank[0_P,z]}=\sum_{i=0}^n (-1)^i f_i(P)\,,\]
where $n=\rank(P)$ and $f_i(P)$ is the number of elements of $P$ of rank $i$ (the ordinary
$f$-vector). Suppose that $K$ is a finite cell complex of
dimension $n-2$ and let $P$ be the poset of non-empty faces of $K$, with a bottom and top
elements adjoined. Then, for $1\leq i\leq n-1$, $f_i(P)$ is the number of faces of $K$ of
dimension $i-1$, and
\[\chi(P)= 1+(-1)^n-\sum_{i=1}^{n-1} (-1)^{i-1} f_i(P)\]
is the Euler characteristic of a sphere of dimension $n-2$ (or $n$)  minus the Euler
characteristic of $K$.

Consider the class of graded posets $P$ for which
\[(\chi-\epsilon)(P)=(\Bar{\zeta}-\zeta^{-1})(P)\,.\]
This includes the graded poset of rank $0$ (trivially) and all posets
$P$ constructed from {\em regular} finite cell complexes as above. This follows
from the above description of $\chi$, together with the fact that the M\"obius
function of $P$ is in this case the reduced Euler characteristic of 
$K$~\cite[Proposition 3.8.8]{Sta86}.
\end{exa}

\begin{exa}\label{Ex:posets} Let $\calP$ be  the
$\F$-vector space with basis the set of all isomorphism classes of finite posets, not
necessarily graded. $\calP$ is a graded connected Hopf algebra where the degree of a poset $P$
is the number of elements of $P$, denoted $\#P$, the product is disjoint union of posets
\[P\cdot Q:=P\sqcup Q\,,\]
the unit element is the poset with no elements, and the coproduct is
\[\Delta(P):=\sum_{I\leq P}I\otimes(P\setminus I)\,.\]
Here, $I\leq P$ indicates that $I$ is a {\em lower ideal} of $P$, i.e., if $x\in P$, $y\in I$ and
$x\leq y$ then $x\in I$. The subsets $I$ and $P\setminus I$ of $P$ are viewed as posets with
the partial order of $P$.

We turn $\calP$ into a combinatorial Hopf algebra by defining $\zeta:\calP\to\F$ by
$\zeta(P)=1$ for every poset $P$. A poset is {\em discrete} if no two elements of $P$ are
comparable. It follows that
\begin{equation}\label{E:mobiusP}
\zeta^{-1}(P)=\begin{cases}(-1)^{\#P} & \text{ if $P$ is
discrete},\\ 0 & \text{ otherwise.}\end{cases}
\end{equation}
In order to see this, note that a subset of $P$ is both discrete and a lower ideal if and only
if it is a subset of $\Min(P)$, the set of minimal elements of $P$. Let $\mu(P)$ be the
functional defined by the right hand side of~\eqref{E:mobiusP}. Then
\[(\mu\zeta)(P)=\sum_{I\leq P}\mu(I)=\sumsub{I\leq P\\I\text{ discrete}}(-1)^{\#I}=
\sum_{I\subseteq\Min(P)}(-1)^{\#I}=\begin{cases}1 & \text{ if
$\Min(P)=\emptyset$},\\ 0 & \text{ otherwise.}\end{cases}\]
Now, any non-empty poset has a minimal element, so $\mu\zeta=\epsilon$ and $\mu=\zeta^{-1}$.

{}From~\eqref{E:mobiusP} we deduce
\[\Bar{\zeta}^{-1}(P)=\begin{cases}1 & \text{ if $P$ is
discrete},\\ 0 & \text{ otherwise.}\end{cases}\]
Hence, the odd character is
\[\nu(P)=\sum_{I\leq P}\Bar{\zeta}^{-1}(I)=\sumsub{I\leq P\\I\text{
discrete}}1=2^{\#\Min(P)}\,.\]
\end{exa}

\begin{exa}\label{Ex:dist-lattice} There is a canonical morphism of combinatorial Hopf algebras
\[J:(\calP,\zeta)\to(\calR,\zeta)\,,\]
where $J(P)$ is the set of lower ideals of $P$, viewed as a poset under inclusion. $J(P)$ is a
graded poset, with bottom element $\emptyset$, top element $P$, and rank function
$\rank(I)=\#I$. The zeta characters are preserved trivially. Hence so are the other canonical
characters.  

It is known that $J(P)$ is in fact a {\em distributive lattice}, and
that any distributive lattice is of this form for a unique poset $P$~\cite[Theorem 3.4.1]{Sta86}. Note
also that if
$P$ is discrete and $\#P=n$ then $J(P)=B_n$, the Boolean poset on $n$ $elements$ (all subsets of
$P$). Since the M\"obius characters are preserved by $J$, we deduce from~\eqref{E:mobiusP} that
the M\"obius function of a distributive lattice $L$ is
\[\mu(L)=\begin{cases}(-1)^{n} & \text{ if $J=B_n$},\\ 0 & \text{ otherwise.}\end{cases}\]
\end{exa}

In Section~\ref{S:terminal}, the canonical characters of the combinatorial Hopf algebra of
quasi-symmetric functions are calculated.

 %%%%%%%%%%%%%%%%%%%%%%%%%%%%%%%%%%%%%%%%%%%%%%%%%%%%%%%%%%%%%%
 \section{The Hopf algebras $\QSym$ and $\NSym$}\label{S:Q-N}

This section  recalls the Hopf algebraic structure of
quasi-symmetric functions and of non-commutative symmetric functions. For more details,
see~\cite{GKal,Ha03,Malv,MR,Sta99}. The reader familiar with these notions
may skip ahead to Section~\ref{S:terminal}.

 A composition $\alpha$ of a
 positive integer $n$, written $\alpha\comp n$, is an ordered list 
$\alpha =(a_1,a_2,\ldots,a_k)$ of positive
 integers such that $a_1+a_2+\cdots+a_k=n$. We let $k(\alpha):=k$ be the number of
parts and $\abs{\alpha}:=n$.
 Compositions of $n$ are in one-to-one correspondence with subsets
 of $\{1,2,\ldots,n-1\}$ via 
\[\alpha\mapsto I(\alpha):=\{a_1,a_1+a_2,\ldots,a_1+\cdots+a_{k-1}\}\,.\]
  For compositions $\alpha,\beta\comp n$, we say that $\alpha$ {\em refines} $\beta$  if
$I(\beta)\subseteq I(\alpha)$, and write $\beta\leq \alpha$.

 The {\em monomial quasi-symmetric function} $M_\alpha$
 indexed by $\alpha=(a_1,\ldots,a_k)$ is
  \begin{equation}\label{E:baseM}
     M_\alpha\ :=\
      \sum_{i_1<i_2<\cdots <i_k} x^{a_1}_{i_1}\cdots x^{a_k}_{i_k}\,.
  \end{equation}
This is an element of the commutative algebra of formal power series in the variables
$\{x_i\}_{i\geq 1}$.  We agree that $M_{()}=1$, where
$()$ denotes the unique composition of $0$ (with no parts).
As $\alpha$ runs over all compositions of $n$, $n\geq 0$, the elements $M_\alpha$ span a
subalgebra $\QSym$ of the algebra of formal power series. It is the algebra of quasi-symmetric
functions. It is in fact a graded subalgebra, the homogeneous component $\QSym_n$ of degree $n$
being spanned by $\{M_\alpha\}_{\alpha\comp n}$. Thus, $\dim\QSym_0=1$ and $\dim
\QSym_n=2^{n-1}$ for $n\geq 1$.

There is a Hopf algebra structure on $\QSym$ with comultiplication
  \begin{equation}\label{E:coprodM}
     \Delta(M_\alpha)\ =\
     \sum_{\alpha =\beta\gamma}M_{\beta }\otimes M_{\gamma},
  \end{equation}
 where $\beta\gamma$ is the concatenation of compositions
 $\beta$ and $\gamma$. The counit is projection onto the (one-dimensional) component of degree
$0$; equivalently, the morphism of algebras that sends all variables $x_i$ to $0$. Thus,
 $\QSym$ is a graded connected bialgebra, and hence a graded Hopf algebra, 
by~\cite{MM65}. An explicit formula for the antipode is given below~\eqref{E:antipode-QSym}.

A second linear basis of $\QSym_n$ (and thus of $\QSym$) is obtained by defining for each
$\alpha\comp n$, 
  \begin{equation}\label{E:baseF}
     F_\alpha\ :=\ \sum_{\alpha\leq \beta} M_\beta\,.
  \end{equation}

\medskip
A partition of $n$ is a composition $\lambda=(l_1,\ldots,l_k)$ of $n$ such that
$l_1\geq\cdots\geq l_k$. Given a composition $\alpha$, let $s(\alpha)$
be the partition obtained by rearranging the parts of $\alpha$ in decreasing order.
The  monomial symmetric functions are the elements
\begin{equation}\label{E:monomial}
m_\lambda\ :=\ \sum_{s(\alpha)=\lambda}M_\alpha\,.
\end{equation}
The subspace spanned by the elements $m_\lambda$, as $\lambda$ runs over all
partitions, is a Hopf subalgebra of $\QSym$. It is the Hopf algebra $\Sym$ of symmetric
functions.

\medskip

 Let $\NSym=\F\langle H_1,H_2,\ldots\rangle$ be the non-commutative
 algebra freely generated by infinitely many variables $\{H_n\}_{n\geq 1}$.
 Define
  \begin{equation}\label{E:coprodH}
     \Delta(H_n)\ :=\ \sum_{i+j=n} H_i\otimes H_j
  \end{equation}
 where $H_0=1$. This turns $\NSym$ into a graded connected Hopf algebra, where
 $\deg(H_n)=n$. A linear basis for the homogeneous component $\NSym_n$ of degree $n$
is $\{H_\alpha\}_{\alpha\comp n}$, where
 \begin{equation}\label{E:prodH}
 H_{(a_1,a_2,\ldots,a_k)}\ :=\ H_{a_1}H_{a_2}\cdots H_{a_k}\,.
\end{equation}

The graded Hopf algebras $\QSym$ and $\NSym$ are dual to each other, in the graded sense.
The identification between $\QSym_n$ and $(\NSym_n)^*$ is
\begin{equation}\label{E:MHduality}
M_\alpha \leftrightarrow H_\alpha^*\,.
\end{equation}

 The ideal ${\mathcal I}$ of $\NSym$ generated by commutators
  is a Hopf ideal and the quotient $\NSym\big/ {\mathcal I}$ is 
 $\Sym$, the  Hopf algebra of symmetric functions.
 The quotient map sends $H_\alpha$ to the complete homogeneous symmetric function
$h_{s(\alpha)}$. 

The Hopf algebra $\Sym$ is self-dual (in the graded sense), the identification being
\begin{equation}\label{E:mhduality}
m_\lambda \leftrightarrow h_\lambda^*\,.
\end{equation}
Therefore, the maps $\NSym\onto\Sym$ and $\Sym\inc\QSym$ are dual to each other.

%%%%%%%%%%%%%%%%%%%%%%%%%%%%%%%%%%%%%%%%%%%%%%%%%%%%%%%%%%%%%%
\section{The terminal object in the category of combinatorial Hopf algebras} \label{S:terminal}

The results of this section do not require any assumptions on the base field $\F$.

We endow $\QSym$ with a canonical character.
Let 
\[\zetaQ:\F[x_1,x_2,\ldots]\to\F\]
 be the morphism of algebras such that
\[\zetaQ(x_1)=1 \text{ and }\zetaQ(x_i)=0 \text{ for all } i\neq 1\,.\]
The map $\zetaQ$ is actually defined on the algebra of power series of finite degree, and hence on $\QSym$.
Note also that all multiplicative functionals on power series which set one
variable to $1$ and the remaining variables to $0$ agree when restricted to $\QSym$, by
quasi-symmetry.

It follows from~\eqref{E:baseM} and~\eqref{E:baseF} that

\begin{equation}\label{E:zeta-QSym}
\zetaQ(M_\alpha)\ =\ \zetaQ(F_\alpha)\ =\ \begin{cases}
       1 & \text{ if $\alpha=(n)$ or $()$,}\\
       0& \text{otherwise.}\end{cases}
\end{equation}

By~\eqref{E:MHduality} we have
\[\zetaQ|_{(\QSym)_n}=H_n\,,\]
as elements of $(\QSym_n)^*=\NSym_n$.
 
 $(\QSym,\zetaQ)$ is a terminal object both as a 
 combinatorial coalgebra and as a combinatorial Hopf algebra.

 \begin{theo}\label{T:univ-qsym}
 For any combinatorial coalgebra (Hopf algebra) $(\H,\zeta)$, there exists a
 unique morphism of combinatorial coalgebras (Hopf algebras) 
\[\Psi:(\H,\zeta)\to(\QSym,\zetaQ)\,.\]
Moreover, $\Psi$ is explicitly given as follows. For $h\in \H_n$,
\begin{equation}\label{E:univ-qsym}
\Psi(h)=\sum_{\alpha\comp n}\zeta_\alpha(h)M_\alpha
\end{equation}
where, for $\alpha=(a_1,\ldots,a_k)$, $\zeta_\alpha$ is the composite
\[\H\map{\Delta^{(k-1)}}\H^{\otimes k}\onto
\H_{a_1}\otimes\cdots\otimes\H_{a_k}\map{\zeta^{\otimes k}}\F\,,\]
where the unlabelled map is the canonical projection onto a homogeneous component.
 \end{theo}
\begin{proof}
Consider the coalgebra case first. As recalled in Section~\ref{S:Q-N}, the graded dual of the
graded coalgebra
$\QSym$ is the free algebra $\NSym$, which has one generator $H_n$ of degree $n$ for every $n$.
Given a graded coalgebra $\H$ and a linear functional $\zeta:\H\to\F$, let (as in
Section~\ref{S:characters})
\[\zeta_n:=\zeta|_{\H_n}:\H_n\to\F\,.\]
Thus $\zeta_n\in (\H_n)^*$, and there is a unique morphism of graded algebras $\Phi:\NSym\to
\H^*$ such that $\Phi(H_n)=\zeta_n$. 

Set $\Psi=\Phi^*:\H\to\QSym$. Thus,
$\Psi$ is a morphism of graded coalgebras and
\[\zetaQ\circ\Psi|_{\H_n}=\zetaQ|_{\QSym_n}\circ\Psi|_{\H_n}=H_n\circ\Psi|_{\H_n}
=\Phi(H_n)=\zeta_n=\zeta|_{\H_n}\,.\]
Thus, $\zetaQ\circ\Psi=\zeta$, proving that $\Psi:\H\to\QSym$ is a morphism of combinatorial
coalgebras. 

Let $\alpha=(a_1,\ldots,a_k)$. By construction, 
$\zeta_\alpha=\zeta_{a_1}\cdots\zeta_{a_k}$ (product in $\H^*$), and by~\eqref{E:prodH},
$H_{\alpha}=H_{a_1}H_{a_2}\cdots H_{a_k}$ (product in $\NSym$). Therefore,
$\Phi(H_\alpha)=\zeta_\alpha$, and hence $\Psi$ is given by~\eqref{E:univ-qsym}.

Uniqueness of $\Psi$ follows from uniqueness of $\Phi$, by duality.

Now suppose that $\H$ is a graded Hopf algebra and $\zeta$ is a morphism of algebras.
It only remains to show that $\Psi$, as constructed above, is a morphism of algebras.
 Consider the following diagrams, where $m$ denotes the multiplication map of $\H$ or of $\QSym$.
  $$
   \begin{picture}(140,58)
    \put(0,45){$\H^{\otimes 2}$}
    \put(60,45){$\H$}
    \put(105,45){$\QSym$}
    \put(60,0){$\F$}
    \put(35,51){\scriptsize $m$}
    \put(85,51){\scriptsize $\Psi$}
    \put(12,25){\scriptsize $\zeta^{\otimes 2}$}
    \put(100,25){\scriptsize$\zetaQ$}
    \put(55,25){\scriptsize$\zeta$}
    \put(27,48){\vector(1,0){27}}
    \put(75,48){\vector(1,0){27}}
    \put(65,40){\vector(0,-1){27}}
    \put(27,40){\vector(1,-1){27}}
    \put(103,40){\vector(-1,-1){27}}
   \end{picture}
 \qquad\quad
   \begin{picture}(200,58)
    \put(0,45){$\H^{\otimes 2}$}
    \put(78,45){$\QSym^{\otimes 2}$}
    \put(168,45){$\QSym$}
    \put(90,0){$\F$}
    \put(45,51){\scriptsize $\Psi^{\otimes 2}$}
    \put(135,51){\scriptsize $m$}
    \put(20,25){\scriptsize $\zeta^{\otimes 2}$}
    \put(155,25){\scriptsize$\zetaQ$}
    \put(77,25){\scriptsize$\zetaQ^{\otimes 2}$}
    \put(29,48){\vector(1,0){43}}
    \put(121,48){\vector(1,0){43}}
    \put(95,40){\vector(0,-1){27}}
    \put(30,40){\vector(2,-1){54}}
    \put(163,40){\vector(-2,-1){54}}
   \end{picture}
  $$
These diagrams commute, by construction of $\Psi$ and since $\zeta$ and $\zetaQ$ are morphisms
of algebras. Also, $\Psi\circ m$ and $m\circ\Psi^{\otimes 2}$ are morphisms of graded coalgebras.
Hence, by the universal property of $\QSym$ as a combinatorial coalgebra applied to
the functional $\zeta^{\otimes 2}$, we have $\Psi\circ m=m\circ\Psi^{\otimes 2}$. Thus, $\Psi$ is
a morphism of algebras.
\end{proof}

\begin{rem} $(\QSym,\zetaQ)$ is also the terminal object in the larger
category of graded coalgebras (Hopf algebras) endowed with a linear (multiplicative) functional,
without any finite-dimensionality assumptions on the homogeneous components. 
A direct proof of this statement (without resorting to the freeness of $\NSym$)
will be given in~\cite{A2}. For most applications of interest, however, 
the present result suffices.
\end{rem}

The Hopf algebra of symmetric functions satisfies the same universal property but among
cocommutative coalgebras. Let $\zetaS:\Sym\to\F$ denote the restriction of
$\zetaQ:\QSym\to\F$. Thus,
\[ \zetaS(m_\lambda)\ =\ \begin{cases}
                    1 & \text{ if $\lambda=(n)$ or $()$,}\\
                0& \text{otherwise.}\end{cases}\]

\begin{theo}\label{T:univ-sym}
For any cocommutative combinatorial coalgebra (Hopf algebra) $(\H,\zeta)$, there exists a
 unique morphism of combinatorial coalgebras (Hopf algebras) 
\[\Psi:(\H,\zeta)\to(\Sym,\zetaS)\,.\]
 \end{theo}
\begin{proof} The above argument can be repeated word by word, as the graded
dual of $\Sym$ ($\Sym$ itself) is the free commutative algebra on the generators
$\{h_n\}_{n\geq 1}$.
\end{proof}

\begin{exa}\label{Ex:flagvector}
 Let $\calR$ be Rota's Hopf algebra
(Example~\ref{Ex:Rota}).  Theorem~\ref{T:univ-qsym}
yields a morphism of graded Hopf algebras $\Psi:\calR\to\QSym$. If $P$ a graded poset of rank
$n$ and $\alpha=(a_1,\ldots,a_k)$ is a composition of $n$ then
$\zeta_\alpha(P)$ is the number of chains 
\[0_P=z_0<z_1<\cdots<z_k=1_P\]
in $P$ such that $\rank[z_{i-1},z_i]=a_i$ for every $i=1,\ldots,k$.
Thus $\zeta_\alpha(P)=f_\alpha(P)$, the {\em flag
$f$-vector} of $P$, and~\eqref{E:univ-qsym} becomes
\[\Psi(P)\ =\ \sum_{\alpha\comp n} f_\alpha(P)
M_\alpha=\sumsub{0_P=z_0<z_1<\cdots<z_k=1_P\\i_1<\cdots<i_k}
\!\!\! x_{i_1}^{\rank[z_0,z_1]}\cdots x_{i_k}^{\rank[z_{k-1},z_k]}\,.\]
 This morphism was introduced by Ehrenborg~\cite{Eh96}.
\end{exa}

\begin{exa}\label{Ex:chromatic} (Chromatic Hopf algebra). Let $\calG$ be the $\F$-vector space
with basis the set of all isomorphism classes of finite (unoriented) graphs. Let $V(G)$ denote
the set of vertices of a graph $G$. For $S\subseteq V(G)$, let $G|_S$ denote the graph with set
of vertices
$S$ and edges those edges of $G$ with both ends in $S$. $\calG$ is a graded connected Hopf
algebra with degree $\abs{G}:=\#V(G)$ (number of vertices), product
$G\cdot H:=G\sqcup H$ (disjoint union of graphs) and coproduct
\[\Delta(G):=\sum_{S\subseteq V(G)}G|_S\otimes G|_{S^c}\,.\]
Note that $\calG$ is cocommutative. This Hopf algebra was considered by Schmitt~\cite{Sch94}.

A graph is discrete if it has no edges. Let $\zeta:\calG\to\F$ be 
\[\zeta(G)=\begin{cases} 1 & \text{ if $G$ is discrete,}\\
0 & \text{ otherwise.} \end{cases}\]
 Thus, $(\calG,\zeta)$ is a cocommutative combinatorial
Hopf algebra. Let $\Psi:\calG\to\Sym$ be the corresponding 
morphism of graded Hopf algebras given by Theorem~\ref{T:univ-sym}. If $G$ is a graph with
$n$ vertices and $\alpha=(a_1,\cdots,a_k)$ is a composition of $n$, then
$\zeta_\alpha(G)$ is the number of ordered decompositions
\[V(G)=S_1\sqcup\cdots\sqcup S_k\]
such that $G|_{S_i}$ is discrete and $\#S_i=a_i$ for every $i=1,\ldots,k$.
A {\em proper coloring} of $G$ is a map
$f:V(G)\to\{1,2,\ldots\}$ such that $f(v)\neq f(w)$ whenever there is an edge in $G$ joining $v$
and $w$. For each ordered decomposition as above and positive integers $i_1<\cdots<i_k$ there is
a proper coloring $f$ given by $f|_{S_j}=i_j$, and conversely.
It follows from~\eqref{E:baseM} and~\eqref{E:univ-qsym} that 
\[\Psi(G)=\sum_{f}\prod_{v\in V(G)}x_{f(v)}\,,\]
the sum over the set of  proper colorings of $G$.  $\Psi(G)$ is Stanley's chromatic
symmetric function~\cite{Sta95}.
\end{exa}

\begin{exa}\label{Ex:Pieri}In~\cite{BMSW,BMSW2,BS99b}, morphisms of Hopf
algebras were constructed from Pieri operations on posets. We explain 
this terminology and the connection to our constructions next.

 Let $P$ be a fixed finite graded poset. Let $\IP$ denote the
subspace of $\calR$ linearly spanned by the (isomorphism classes of) intervals of $P$.
Then $\IP$ is a  graded subcoalgebra of $\calR$ (in the terminology of Joni and Rota,
this is the {\em reduced incidence coalgebra} of $P$). Let $\F P$ be the vector space with
basis the underlying set of $P$. This carries a natural left graded comodule
structure over $\IP$,
\[\chi:\F P\to \IP\otimes \F P,\ \ \ \ \chi(y)= \sum_{y\leq z}[y,z]\otimes z\,,\]
where the degree of $z\in P$ is $\rank[0_p,z]$.

Now suppose that a linear functional $\zeta:\IP\to\F$ is given. 
By Theorem~\ref{T:univ-qsym}, there is a corresponding morphism of graded coalgebras
$\Psi:\IP\to\QSym$. This allows us to view $\F P$ as a left graded comodule over
$\QSym$, via
\[\F P\map{\chi} \IP\otimes \F P\map{\Psi\otimes\id} \QSym\otimes\F P,\ \ \ \ 
y\mapsto \sum_{y\leq z}\Psi([y,z])\otimes z\,.\]
Dualizing, we obtain a right graded module structure on $\F P$ over $\NSym$, which
is determined by
\[y\cdot H_n = \sumsub{y\leq z\\\rank[y,z]=n} \zeta([y,z]) \cdot z\,.\]
In the terminology of~\cite{BMSW}, this right module structure constitutes a {\em Pieri
operation} on $P$. Note that specifying such a structure is equivalent to
giving the linear functional $\zeta$. The morphism of Hopf algebras
constructed in~\cite{BMSW} is $\Psi$. Thus, the construction of Hopf algebra morphisms from
Pieri operations is another instance of the general construction of Theorem~\ref{T:univ-qsym}.
\end{exa}

Many other examples of a combinatorial nature will be presented in~\cite{A2}. Additional 
applications of the universal property are discussed below
(Examples~\ref{Ex:powerzeta}--\ref{Ex:theta}).

Consider now the case when $\H=\QSym$ and $\zeta$ is an arbitrary character on $\QSym$. Let
$\alpha$ and $\beta$ be two compositions of $n$. It follows from~\eqref{E:coprodM} that
\[\zeta_\alpha(M_\beta)=\begin{cases}
\zeta(M_{\beta_1})\cdots\zeta(M_{\beta_h}) & \text{ if }\alpha\leq\beta\,,\\
0 & \text{ if not}
\end{cases}\]
where, for $\alpha\leq\beta$, the compositions
$\beta_1,\ldots,\beta_h$ are defined by $\beta=\beta_1\cdots\beta_h$ (concatenation of
compositions) and
$\alpha=(\abs{\beta_1},\ldots,\abs{\beta_h})$. Formula~\eqref{E:univ-qsym} becomes
\begin{equation}\label{E:qsym-qsym}
\Psi(M_\beta)\ = \ \sum_{\beta=\beta_1\cdots\beta_h}\zeta(M_{\beta_1})\cdots\zeta(M_{\beta_h})
M_{(\abs{\beta_1},\ldots,\abs{\beta_h})}\,.
\end{equation}
In particular, $\Psi$ is triangular with respect to refinement of compositions, and $\Psi$ is
an isomorphism if and only if $\Psi(M_{(n)})\neq 0$ for every $n$.

\begin{exa} \label{Ex:powerzeta} Let $m\in\ZZ$ be an integer and $\Psi_m:\QSym\to\QSym$ be the
unique morphism of Hopf algebras such that
$\zetaQ\circ\Psi_m=\zetaQ^m$, the $m$-th convolution power of $\zetaQ$. It follows
from~\eqref{E:coprodM} that for any composition $\beta$ and non-negative integer $m$,
\begin{equation}\label{E:powerzeta}
\zetaQ^m(M_\beta)=\binom{m}{k(\beta)}\,.
\end{equation}
(For general reasons, this formula remains valid for negative integers $m$ as well.)

Then, by~\eqref{E:qsym-qsym},
\begin{equation}\label{E:powerpsi}
\Psi_m(M_\beta)=\sum_{\beta=\beta_1\cdots\beta_h}\binom{m}{k(\beta_1)}\cdots
\binom{m}{k(\beta_h)}M_{(\abs{\beta_1},\ldots,\abs{\beta_h})}\,.
\end{equation}
In particular, $\Psi_1=\id$ (as had to be the case by uniqueness) and
\[\Psi_{-1}(M_\beta)=(-1)^{k(\beta)}\sum_{\alpha\leq\beta}M_\alpha\,.\]
\end{exa}

\begin{exa}\label{Ex:qsym-antipode} Recall that the {\em cooposite} coalgebra $C^{\cop}$ of a
coalgebra $C$ is obtained by reversing the order of the tensor factors in the coproduct of
$C$. Let us apply the universal property to
$(\H,\zeta)=(\QSym^{\cop},\zetaQ^{-1})$.  The resulting morphism of Hopf algebras
$\QSym^{\cop}\to\QSym$ is the antipode $S$ of
$\QSym$, since $S:\QSym\to\QSym$ is an antimorphism of coalgebras and $\zetaQ\circ S=\zetaQ^{-1}$
(see Section~\ref{S:characters}). Theorem~\ref{T:univ-qsym} will thus furnish an explicit
formula for the antipode of $\QSym$. {}From~\eqref{E:powerzeta}, and since
$\binom{-1}{k}=(-1)^k$, we get
\begin{equation}\label{E:mu-QSym}
\zetaQ^{-1}(M_\beta)=(-1)^{k(\beta)}\,.
\end{equation}
(This may also be obtained directly from the recursion $\zetaQ^{-1}\zetaQ=\epsilon$.)

Then, since $\Delta^{\cop}(M_{(a_1,\ldots,a_k)})=\sum_{i=0}^k M_{(a_{i+1},\ldots,a_1)}\otimes
M_{(a_1,\ldots,a_i)}$, formula~\eqref{E:univ-qsym} yields
\begin{equation}\label{E:antipode-QSym}
S(M_\beta)\ =\ (-1)^{k(\beta)}\sum_{\alpha\leq\widetilde{\beta}}M_\alpha\,,
\end{equation}
where if $\beta=(b_1,b_2,\ldots,b_k)$ then 
$\widetilde{\beta}=(b_k,\ldots,b_2,b_1)$. This formula appears in~\cite{Eh96}
and~\cite{Malv}.
\end{exa}

\begin{exa}\label{Ex:theta} 
Consider the odd character $\nuQ=\BzetaQ^{-1}\zetaQ$
(Example~\ref{Ex:nu-chi}).  It follows from~\eqref{E:mu-QSym} that
$\BzetaQ^{-1}(M_\beta)=(-1)^{\abs{\beta}+k(\beta)}$. Let $\beta=(b_1,\ldots,b_k)$. In
\[\nuQ(M_\beta)=\sum_{i=0}^k \BzetaQ^{-1}(M_{(b_1,\ldots,b_i)})\zetaQ(M_{(b_{i+1},\ldots,b_k)})\]
only the terms corresponding to $i=k-1$ and $i=k$ are non zero. These cancel each other if $b_k$
is even and are equal to $(-1)^{\abs{\beta}+k(\beta)}$ if $b_k$ is odd. Thus,
\begin{equation}\label{E:nu-QSym}
\nuQ(M_{\beta})=\begin{cases}
1 & \text{ if $\alpha=()$,}\\
2\cdot(-1)^{\abs{\beta}+k(\beta)} & \text{ if the last part of $\beta$ is odd,}\\
0 & \text{ otherwise.}
\end{cases}
\end{equation}

Let $\Theta:\QSym\to\QSym$ be the unique morphism of Hopf algebras such that
$\zetaQ\circ\Theta=\nuQ$. By~\eqref{E:qsym-qsym},
\[\Theta(M_{\beta})=\sum_{\beta=\beta_1\cdots\beta_h}\nuQ(M_{\beta_1})\cdots\nuQ(M_{\beta_h})
M_{(\abs{\beta_1},\ldots,\abs{\beta_h})}\,.\]
Let $\alpha=(\abs{\beta_1},\ldots,\abs{\beta_h})$. For $M_\alpha$ to appear in the above sum,
the last part of each $\beta_i$ must be odd. In this case, $M_\alpha$ appears with coefficient
\[2^{\sum_{i=1}^h 1}\cdot
(-1)^{\sum_{i=1}^h\abs{\beta_i}+k(\beta_i)}=2^{k(\alpha)}\cdot(-1)^{\abs{\beta}+k(\beta)}\,.\]
Thus, 
\[\Theta(M_\beta)=(-1)^{\abs{\beta}+k(\beta)}\cdot\sum_{\alpha}2^{k(\alpha)}M_\alpha\,,\]
the sum over those $\alpha\leq\beta$ obtained by adding the entries within consecutive segments
$\beta_1,\ldots,\beta_h$ of $\beta$ with the property that the last part of each $\beta_i$ is
odd. This set of compositions is empty if the last part of $\beta$ is even.
Otherwise, it is a lower interval whose maximum element is the composition obtained by adding the
entries within each maximal segment of $\beta$ of the form (even,\,even,\ldots,odd). Let
$\odd(\beta)$ denote this element. For instance,
$\odd(1,2,2,1,3,2,3)=(1,5,3,5)$. It follows that
\begin{equation}\label{E:theta}
\Theta(M_\beta)=\begin{cases}
1 & \text{ if $\beta=()$,}\\
(-1)^{\abs{\beta}+k(\beta)}\cdot\sum_{\alpha\leq\odd(\beta)}2^{k(\alpha)}M_\alpha
& \text{ if the last part of $\beta$ is odd,}\\
0 & \text{ otherwise.}
\end{cases}
\end{equation}
This shows that $\Theta$ is the morphism of Hopf algebras introduced by
Stembridge~\cite[Theorem 3.1]{Ste97}. This map was originally defined in the
$F$-basis of $\QSym$. {}From here, Hsiao derived an expression for $\Theta$ in the
$M$-basis~\cite[Theorem 2.4]{Hsiao}, which agrees with~\eqref{E:theta}.
The original expression may also be easily obtained through the universal property of $\QSym$,
starting from the fact that
\[\nuQ(F_\alpha)\ =\ \begin{cases}
                   1 & \text{ if $\alpha=()$,}\\
              2& \text{ if $\alpha=(1,1,\ldots,1,k)$,}\\
  0& \text{ otherwise.}\end{cases} \]

\end{exa}

Other characters on $\QSym$ and the corresponding morphisms are discussed in
Section~\ref{S:even-odd-QSym}.

%%%%%%%%%%%%%%%%%%%%%%%%%%%%%%%%%%%%%%%%%%%%%%%%%%%%%%%%%%%%%%%%%
\section{Even and odd subalgebras of a combinatorial Hopf algebra}\label{S:even-odd}

Let $\H$ be a graded connected Hopf algebra. Given a linear functional $\varphi:\H\to\F$
we let $\varphi_n:=\varphi|_{\H_n}$. This is an element of degree $n$ of the graded dual
$\H^*$. 

\begin{defi}\label{D:S-I}
 Given characters $\varphi,\psi:\H\to\F$, we let
$S(\varphi,\psi)$ denote the largest graded subcoalgebra of $\H$ such that 
\[\forall\ h\in S(\varphi,\psi),\ \ \varphi(h)=\psi(h)\,,\]
and we let $I(\varphi,\psi)$ denote the ideal of $\H^*$ generated by
\[\varphi_n-\psi_n,\ \ n\geq 0\,.\]
\end{defi}
\begin{rem} By definition, $S(\varphi,\psi)$ is the largest
{\em graded} subcoalgebra contained in $\ker(\varphi-\psi)$ (a subspace of codimension $\leq 1$).
It is easy to see that $S(\varphi,\psi)$ may also be defined as the largest subcoalgebra (graded
or not) contained in the graded subspace $\oplus_{n\geq 0}\ker(\varphi_n-\psi_n)$ (this is
the approach taken in~\cite{A}). We will not make use of this fact.
\end{rem}

The basic properties of these objects follow. 

\begin{theo}\label{T:S-I}\hfill
\begin{itemize}
\item[(a)] $S(\varphi,\psi)=\{h\in\H\ \mid\ f(h)=0\ \forall\ f\in I(\varphi,\psi)\}$.
\item[(b)] $I(\varphi,\psi)$ is a graded Hopf ideal of $\H^*$.
\item[(c)] $S(\varphi,\psi)$ is a graded Hopf subalgebra of $\H$.
\item[(d)] A homogeneous element $h\in\H$ belongs to $S(\varphi,\psi)$ if and only if
\begin{equation}\label{E:S-explicit}
\bigl(\id\otimes(\varphi-\psi)\otimes\id\bigr)\circ\Delta^{(2)}(h)=0\,.
\end{equation}
\end{itemize}
\end{theo}
\begin{proof} 
Write $I(\varphi,\psi)^\perp:=\{h\in\H\ \mid\ f(h)=0\ \forall\ f\in I(\varphi,\psi)\}$.
Since $I(\varphi,\psi)$ is generated by homogeneous elements, it is a graded ideal of $\H^*$. 
Hence $I(\varphi,\psi)^\perp$ is a graded subcoalgebra of $\H$. 

If $h\in
I(\varphi,\psi)^\perp$ then 
\[\varphi(h)=\sum_n\varphi_n(h)=\sum_n\psi_n(h)= \psi(h)\,,\] 
since $\varphi_n-\psi_n\in I(\varphi,\psi)$. Therefore, 
$I(\varphi,\psi)^\perp\subseteq S(\varphi,\psi)$. 

Let $C$ be a graded subcoalgebra of $\H$
such that $\varphi(h)=\psi(h)$ for all $h\in C$. Let $h\in C$. Write $h=\sum_n h_n$. Since $C$
is graded, each $h_n\in C$. Hence, 
\[(\varphi_n-\psi_n)(h)=(\varphi_n-\psi_n)(h_n)=(\varphi-\psi)(h_n)=0\,.\]
Thus, $(\varphi_n-\psi_n)(C)=0$ for all $n\geq 0$. Since $C$ is a subcoalgebra and
$I(\varphi,\psi)$ is the ideal generated by $\varphi_n-\psi_n$, $n\geq 0$, it follows that
$f(C)=0$ for all $f\in I(\varphi,\psi)$. Thus, $C\subseteq I(\varphi,\psi)^\perp$ and
hence
$S(\varphi,\psi)\subseteq I(\varphi,\psi)^\perp$. The proof of (a) is complete.

If $C$ and $D$ are graded subcoalgebras of $\H$, then so is their product $C\cdot D$. Also, since
$\varphi$ and $\psi$ are characters we have
\[(\varphi-\psi)(xy)=(\varphi-\psi)(x)\varphi(y)+\psi(x)(\varphi-\psi)(y)\,.\]
Hence, if $(\varphi-\psi)(C)=(\varphi-\psi)(D)=0$, then $(\varphi-\psi)(C\cdot D)=0$.
It follows that $S(\varphi,\psi)\cdot S(\varphi,\psi)\subseteq S(\varphi,\psi)$.
Also, $\H_0=\F\cdot 1$ is a graded subcoalgebra of $\H$ and $\varphi(1)=1=\psi(1)$,
so it is contained in $S(\varphi,\psi)$. Thus, $S(\varphi,\psi)$ is a subalgebra of
$\H$ and by (a) $I(\varphi,\psi)$ is a coideal of $\H^*$. This completes the proofs of (b) and
(c).

Let $h$ be a homogeneous element satisfying~\eqref{E:S-explicit}. Then
\[\sum_{n\geq 0}\bigl(\id\otimes(\varphi_n-\psi_n)\otimes\id\bigr)\circ\Delta^{(2)}(h)=0\,.\]
Each term in this sum is homogeneous and their degrees are distinct. Therefore, for each $n\geq
0$,
\[\bigl(\id\otimes(\varphi_n-\psi_n)\otimes\id\bigr)\circ\Delta^{(2)}(h)=0\,.\]
Let $\lambda,\rho$ be arbitrary elements of $\H^*$. We deduce that for each $n\geq
0$,
\[\bigl(\lambda(\varphi_n-\psi_n)\rho\bigr)(h)=\bigl(\lambda\otimes(\varphi_n-\psi_n)\otimes\rho\bigr)\circ\Delta^{(2)}(h)=0\,.\]
It follows that for any $f\in I(\varphi,\psi)$, $f(h)=0$. Thus, by (a), $h\in S(\varphi,\psi)$.

Conversely, let $C$ be a graded subcoalgebra of $\H$ such that $(\varphi-\psi)(C)=0$.
Then
\[\bigl(\id\otimes(\varphi-\psi)\otimes\id\bigr)\circ\Delta^{(2)}(C)\subseteq
C\otimes(\varphi-\psi)(C)\otimes C=0\,.\]
Thus, every homogeneous element of $C$ satisfies~\eqref{E:S-explicit}. This applies in particular
to $C=S(\varphi,\psi)$, which completes the proof of (d).
\end{proof}

\begin{coro}\label{C:S-I-duality} There is an isomorphism of graded Hopf algebras
\[S(\varphi,\psi)^*\cong \H/I(\varphi,\psi)\,.\]
\end{coro}
\begin{proof} This follows at once from Theorem~\ref{T:S-I}.
\end{proof}

\begin{prop}\label{P:S-I-invariance}
Let $\varphi,\varphi',\psi,\psi'$ be characters on $\H$. Suppose that either
\[\psi^{-1}\varphi=(\psi')^{-1}\varphi' \text{ \ or \ }
\varphi\psi^{-1}=\varphi'(\psi')^{-1}\,.\]
Then 
\[S(\varphi,\psi)=S(\varphi',\psi') \text{ \ and \ }
I(\varphi,\psi)=I(\varphi',\psi')\,.\]
\end{prop}
\begin{proof} Suppose $\psi^{-1}\varphi=(\psi')^{-1}\varphi'$. Then
\[\psi^{-1}(\varphi-\psi)=\psi^{-1}\varphi-\epsilon=(\psi')^{-1}\varphi'-\epsilon=
(\psi')^{-1}(\varphi'-\psi')\,.\]
Passing to homogeneous components we see that $\varphi_n-\psi_n$ belongs to the (left) ideal
generated by $\varphi'_n-\psi'_n$, $n\geq 0$.
Thus $I(\varphi,\psi)=I(\varphi',\psi')$ and, by Theorem~\ref{T:S-I},
$S(\varphi,\psi)=S(\varphi',\psi')$. The other case follows similarly.
\end{proof}

Note also that $S(\varphi,\psi)=S(\psi,\varphi)$ and $I(\varphi,\psi)=I(\psi,\varphi)$
(symmetry).

\begin{prop}\label{P:S-I-functoriality} Let $\alpha:\H'\to\H$ is a morphism of graded
connected Hopf algebras. Let $\varphi, \psi$ be characters of $\H$ and 
$\varphi':=\varphi\circ\alpha$,
$\psi':=\psi\circ\alpha$. Then,
\begin{itemize}
\item[(a)] $\alpha\bigl(S(\varphi',\psi')\bigr)\subseteq S(\varphi,\psi)$;
\item[(b)] $I(\varphi',\psi')$ is the ideal of $(\H')^*$
generated by $\alpha^*\bigl(I(\varphi,\psi)\bigr)$.
\end{itemize}
Moreover, if $\alpha$ is injective then 
\[S(\varphi',\psi')=\alpha^{-1}\bigr(S(\varphi,\psi)\bigr) \text{ \ and \ }
I(\varphi',\psi')=\alpha^*\bigl(I(\varphi,\psi)\bigr)\,.\]
\end{prop}
\begin{proof} By hypothesis, $\alpha\bigl(\ker(\varphi'-\psi')\bigr)\subseteq\ker(\varphi-\psi)$.
Since the image under $\alpha$ of a graded subcoalgebra is another graded subcoalgebra, 
$\alpha\bigl(S(\varphi',\psi')\bigr)$ is a graded subcoalgebra
contained in $\ker(\varphi-\psi)$, and hence also in $S(\varphi,\psi)$. This proves (a)

Since $\alpha^*(\varphi_n-\psi_n)=\varphi'_n-\psi'_n$, $\alpha^*$ carries the generators of the
ideal $I(\varphi,\psi)$ of $\H^*$ to the generators of the ideal
$I(\varphi',\psi')$ of $(\H')^*$. This implies (b).

If $\alpha$ is injective then the preimage under $\alpha$ of a subcoalgebra of $\H$ is
a subcoalgebra of $\H'$. Therefore, $\alpha^{-1}\bigr(S(\varphi,\psi)\bigr)$ is a graded
subcoalgebra of $\H'$ contained in
$\alpha^{-1}\bigr(\ker(\varphi-\psi)\bigr)=\ker(\varphi'-\psi')$, and hence also in 
$S(\varphi',\psi')$. 

Finally, in this situation $\alpha^*$ is surjective, so the image of an ideal of $\H^*$ under
$\alpha^*$ is an ideal of $(\H')^*$. Hence, by (b),
$I(\varphi',\psi')=\alpha^*\bigl(I(\varphi,\psi)\bigr)$.
\end{proof}

We now specialize these constructions. 

\begin{defi}\label{D:evenodd}
Let $(\H,\zeta)$ be a combinatorial Hopf algebra. The {\em even} and {\em odd subalgebras}
of $\H$ are respectively
\[S_+(\H,\zeta):=S(\Bar{\zeta},\zeta) \text{ \ and \ }
S_-(\H,\zeta):=S(\Bar{\zeta},\zeta^{-1})\,.\]
The {\em even} and {\em odd ideals}
of $\H^*$ are respectively
\[I_+(\H,\zeta):=I(\Bar{\zeta},\zeta) \text{ \ and \ }
I_-(\H,\zeta):=I(\Bar{\zeta},\zeta^{-1})\,.\]
\end{defi}

Specializing the previous results we obtain the basic properties of even and odd subalgebras
and ideals.

\begin{prop}\label{P:even-odd-basic} Let $(\H,\zeta)$ be a combinatorial Hopf algebra.
\begin{itemize}
\item[(a)] $S_+(\H,\zeta)$ and $S_-(\H,\zeta)$ are graded Hopf subalgebras of $\H$,
$I_+(\H,\zeta)$ and $I_-(\H,\zeta)$ are graded Hopf ideals of $\H^*$, and
there are isomorphisms of graded Hopf algebras
\[S_{\pm}(\H,\zeta)^*\cong \H^*/I_{\pm}(\H,\zeta)\,.\]
\item[(b)] $S_+(\H,\zeta)$ and $S_-(\H,\zeta)$ are the largest subcoalgebras  with the
properties that
\[\zeta(h)=(-1)^n\zeta(h) \text{ \ (respectively \ 
$\zeta^{-1}(h)=(-1)^n\zeta(h)$)}
\]
 for every $h\in S_+(\H,\zeta)$ (respectively $S_-(\H,\zeta)$) of degree $n$.
\item[(c)] $\zeta$ is even if and only if $S_+(\H,\zeta)=\H$ and $\zeta$ is odd if
and only if $S_-(\H,\zeta)=\H$. 
\item[(d)] Let $\nu=\Bar{\zeta}^{-1}\zeta$ and $\chi=\Bar{\zeta}\zeta$
be the canonical characters of Example~\ref{Ex:nu-chi}. We have
\[S_+(\H,\zeta)=S(\nu,\epsilon) \text{ \ and \ }S_-(\H,\zeta)=S(\chi,\epsilon)\]
and similarly for $I_{\pm}$.
\item[(e)] Let $\alpha:(\H',\zeta')\to (\H,\zeta)$ be a morphism of
combinatorial Hopf algebras. Then
\[\alpha\bigl(S_{\pm}(\H',\zeta')\bigr)\subseteq S_{\pm}(\H,\zeta) \]
and $I_{\pm}(\H',\zeta')$ is the ideal of $(\H')^*$ generated by
$\alpha^*\bigl(I_{\pm}(\H,\zeta)\bigr)$.
\item[(f)] If in addition $\alpha$ is injective then
\[S_{\pm}(\H',\zeta')=\alpha^{-1}\bigl(S_{\pm}(\H,\zeta)\bigr) \text{ \ and \ }
I_{\pm}(\H',\zeta')=\alpha^*\bigl(I_{\pm}(\H,\zeta)\bigr)\,.\]
\end{itemize}
\end{prop}
\begin{proof} These assertions are special cases of Theorem~\ref{T:S-I},
Corollary~\ref{C:S-I-duality}, and Propositions~\ref{P:S-I-invariance}
and~\ref{P:S-I-functoriality}.
\end{proof}

Even and odd characters of $\H$ are related to the even and odd subalgebras of the dual $\H^*$,
as we now explain. Let $\H$ be a graded connected Hopf algebra, $\zeta:\H\to\F$ a character on $\H$ and
$\eta:\H^*\to\F$ a character on its dual. We do not assume any compatibility between $\zeta$ and
$\eta$. Note that each component $\zeta_n=\zeta|_{\H_n}$ is a homogeneous element of $\H^*$.

\begin{prop}\label{P:duality} In the above situation, if $\zeta$ is even (respectively, odd) then
each component $\zeta_n$ belongs to  $S_+(\H^*,\eta)$ (respectively, $S_-(\H^*,\eta)$).
\end{prop}
\begin{proof}
Observe first that
\[\Bar{\eta}(\zeta_n)=(-1)^n\eta(\zeta_n)=\eta(\Bar{\zeta}_n)\]
and
\[\eta^{-1}(\zeta_n)=(\eta\circ S_{\H^*})(\zeta_n)=\eta(\zeta_n\circ
S_\H)=\eta(\zeta^{-1}_n)\,.\]
Recall that characters of $\H$ are group-like elements of $\H^\circ$. Hence,
\[\Delta_{\H^*}(\zeta_n)=\sum_{i+j=n}\zeta_i\otimes\zeta_j\,.\]
Suppose now that $\zeta$ is odd. We have
\begin{align*}
\bigl(\id\otimes(\Bar{\eta}-\eta^{-1})\otimes\id\bigr)\circ\Delta^{(2)}_{\H^*}(\zeta_n) &=
\sum_{i+j+h=n}\zeta_i\otimes(\Bar{\eta}-\eta^{-1})(\zeta_j)\otimes\zeta_h\\
&=\sum_{i+j+h=n}\zeta_i\otimes\bigl(\eta(\Bar{\zeta}_j)-\eta(\zeta^{-1}_j)\bigr)\otimes\zeta_h=0\,,
\end{align*}
since $\Bar{\zeta}_j=\zeta^{-1}_j$ for every $j$. Thus, by item (d) of Theorem~\ref{T:S-I},
$\zeta_n\in S_-(\H^*,\eta)$. The even case is similar.
\end{proof}

\subsection*{Generalized Dehn-Sommerville relations} Let $(\H,\zeta)$ be a combinatorial Hopf
algebra. Consider its odd subalgebra
\[S_-(\H,\zeta)=S(\Bar{\zeta},\zeta^{-1})=S(\chi,\epsilon)\,.\]
Let $h\in\H$ be a homogeneous element. According to~\eqref{E:S-explicit}, the following
conditions are equivalent:
\begin{gather}
\notag h\in S_-(\H,\zeta)\\
\label{E:DS} \bigl(\id\otimes(\Bar{\zeta}-\zeta^{-1})\otimes\id\bigr)\circ\Delta^{(2)}(h)=0\\
\label{E:BB} \bigl(\id\otimes(\chi-\epsilon)\otimes\id\bigr)\circ\Delta^{(2)}(h)=0\,.
\end{gather}

We refer to either set of equations~\eqref{E:DS} or~\eqref{E:BB} as the {\em generalized
Dehn-Sommerville relations} for the combinatorial Hopf algebra $(\H,\zeta)$.
As they both define the same subalgebra of $\H$, they have the same solutions.

\begin{exa}\label{Ex:BB} We describe the generalized Dehn-Sommerville relations for the
combinatorial Hopf algebra $(\QSym,\zetaQ)$. First, from~\eqref{E:mu-QSym} we get
\begin{equation}\label{E:xi-QSym}
(\BzetaQ-\zetaQ^{-1})(M_\beta)=\begin{cases}0 & \text{ if }\beta=()\\
(-1)^{n}+1 & \text{ if }\beta=(n)\\ 
-(-1)^{k(\beta)} & \text{ otherwise.}\end{cases}
\end{equation}
Given an element $h=\sum_{\gamma\comp n}f_\gamma(h) M_\gamma\in\QSym_n$, a composition 
$\alpha=(a_1,\ldots,a_k)\comp n$, and a number $i\in\{1,\ldots,k\}$, set
\[E_{\alpha,i}(h):=(-1)^{a_i}f_{\alpha}(h)-\sum_{\beta\comp
a_i}(-1)^{k(\beta)}f_{\alpha\circ_i\beta}(h)\,.\]
Here, $\alpha\circ_i\beta$ denotes the composition of $n$ obtained by replacing the $i$th part of
$\alpha$ by the parts of $\beta$: if $\beta=(b_1,\ldots,b_h)$,
\[\alpha\circ_i\beta:=(a_1,\ldots,a_{i-1},b_1,\ldots,b_h,a_{i+1},\ldots,a_k)\,.\]
Also, write
\[\alpha_i:=(a_1,\ldots,a_{i-1}) \text{ \ and \ } \alpha^i:=(a_{i+1},\ldots,a_k) \,.\]
It follows from~\eqref{E:xi-QSym} and~\eqref{E:coprodM} that
\[ \bigl(\id\otimes(\BzetaQ-\zetaQ^{-1})\otimes\id\bigr)\circ\Delta^{(2)}(h)
=\sum_{\alpha\comp n}\sum_{i=1}^{k(\alpha)} E_{\alpha,i}(h)\cdot M_{\alpha_i}\otimes
M_{\alpha^i}\,.\]
Since this is a sum of linealy independent tensors, the generalized Dehn-Sommerville 
relations~\eqref{E:DS} for the element $h$ are
\[E_{\alpha,i}(h)\ =\ 0 \text{ \ for each \ }\alpha\comp n \text{ and }
i\in\{1,\ldots,k(\alpha)\}\,,\] or explicitly
\[(-1)^{a_i}f_{\alpha}(h)\ =\
\sum_{\beta\comp a_i}(-1)^{k(\beta)}f_{\alpha\circ_i\beta}(h)\,.\]
Thus, an element $h\in\QSym$ of degree $n$ belongs to $S_-(\QSym,\zetaQ)$ if and only if these
relations are satisfied. This form of the generalized Dehn-Sommerville relations was introduced
in~\cite[Section 7]{A}.

Now consider~\eqref{E:BB}, the second form of the generalized Dehn-Sommerville relations.
{}From~\eqref{E:zeta-QSym} we get
\begin{equation}\label{E:chi-QSym}
(\chiQ-\epsilonQ)(M_\beta)=\begin{cases}  1+(-1)^{n}& \text{ if }\beta=(n)\\ 
(-1)^{i} & \text{ if }\beta=(i,n-i)\\
0  & \text{ otherwise.}\end{cases}
\end{equation}
Therefore,
\[ \bigl(\id\otimes(\chiQ-\epsilonQ)\otimes\id\bigr)\circ\Delta^{(2)}(h)
=\sum_{\alpha\comp n}\sum_{i=1}^{k(\alpha)} B_{\alpha,i}(h)\cdot M_{\alpha_i}\otimes
M_{\alpha^i}\,,\]
with
\[B_{\alpha,i}(h)\ =\
\sum_{j=0}^{a_i}(-1)^jf_{(a_1,\ldots,a_{i-1},j,a_i-j,a_{i+1},\ldots,a_k)}(h)\,,
\] where it is understood that zero parts are omitted. 
Thus, relations~\eqref{E:BB} for the element $h$ are
\[\sum_{j=0}^{a_i}(-1)^jf_{(a_1,\ldots,a_{i-1},j,a_i-j,a_{i+1},\ldots,a_k)}(h)\ =\ 0 \text{ \
for each \ }(a_1,\ldots,a_k)\comp n,\ i\in\{1,\ldots,k\}\,.\] 
These are precisely the original generalized Dehn-Sommerville relations as introduced by Bayer
and Billera~\cite[Theorem 2.1]{BaBi}. 
\end{exa}

Consider the combinatorial Hopf algebra $(\calR,\zeta)$ of
Example~\ref{Ex:Rota}. 
Let $\calE$ be the subspace of $\calR$ spanned by all eulerian posets. Recall that a
finite graded poset $P$ is  {\em eulerian} if its M\"obius function
satisfies $\mu([x,y])=(-1)^{\rank[x,y]}$ for every
$x<y$ in $P$. In particular, intervals of an eulerian poset and cartesian products of eulerian
posets are again eulerian, so $\calE$ is a Hopf subalgebra of $\calR$. 

According to Example~\ref{Ex:Rota}, for any $P\in\calE$ we have
\[\zeta^{-1}(P)=\Bar{\zeta}(P)\,.\]
Therefore, $\calE$ must be contained in the odd subalgebra of $\calR$:
\[\calE\subseteq S_-(\calR,\zeta)\,.\]
Let $\Psi:(\calR,\zeta)\to(\QSym,\zetaQ)$ be the unique morphism of combinatorial Hopf algebras,
as in Example~\ref{Ex:flagvector}:
\[\Psi(P)\ =\ \sum_{\alpha\comp n} f_\alpha(P)M_\alpha\,,\]
where $f_\alpha(P)$ is the flag vector of $P$.
 By item (e) of Proposition~\ref{P:even-odd-basic},
\[\Psi(\calE)\subseteq \Psi\bigl(S_-(\calR,\zeta)\bigr)\subseteq S_-(\QSym,\zetaQ)\,.\]
We thus recover the important result of Bayer and Billera~\cite[Theorem 2.1]{BaBi} that the flag
vector of any eulerian poset must satisfy the generalized Dehn-Sommerville equations of
Example~\ref{Ex:BB}.

 %%%%%%%%%%%%%%%%%%%%%%%%%%%%%%%%%%%%%%%%%%%%%%%%%%%%%%%%%%%%%%%%%
 \section{Even and odd subalgebras of $\QSym$}\label{S:even-odd-QSym}

Below we describe the even and odd subalgebras of the terminal combinatorial Hopf algebra 
$(\QSym,\zetaQ)$ in explicit terms. We continue to assume that $\ch\F\neq 2$. 
We find that the odd subalgebra is the Hopf algebra
introduced by Stembridge in~\cite{Ste97}. 

Let $\Pi_{+}$ and $\Pi_{-}$ denote the even and odd subalgebras of $(\QSym,\zetaQ)$.

\begin{prop}\label{P:E-Pi-functoriality} Let $(\H,\zeta)$ be a combinatorial 
Hopf algebra and
$\Psi:(\H,\zeta)\to(\QSym,\zetaQ)$ the unique morphism. Then
\[\Psi\bigl(S_+(\H,\zeta)\bigr)\subseteq \Pi_{+} \text{ \ and \ }
\Psi\bigl(S_-(\H,\zeta)\bigr)\subseteq \Pi_{-}\,.\]
\end{prop}
\begin{proof}
This is a particular case of Proposition~\ref{P:even-odd-basic}, item (e).
\end{proof}

\begin{coro}\label{C:E-Pi-univ} Let $(\H,\zeta)$ be a combinatorial Hopf algebra and
consider the unique morphism $\Psi:(\H,\zeta)\to(\QSym,\zetaQ)$. If $\zeta$
is even then $\im(\Psi)\subseteq \Pi_{+}$ and if $\zeta$ is odd then
$\im(\Psi)\subseteq \Pi_{-}$.
In addition, $(\Pi_{+},\zetaQ)$ and $(\Pi_{-},\zetaQ)$ are the terminal objects in the categories of
even and odd combinatorial Hopf algebras, respectively.
\end{coro}
\begin{proof} If $\zeta$ is even then, by item (c) in Proposition~\ref{P:even-odd-basic},
$S_+(\H,\zeta)=\H$, so the conclusion follows from Proposition~\ref{P:E-Pi-functoriality}. The
rest is similar.
\end{proof}

 A composition $\alpha=(a_1,\ldots,a_k)$ is
said to be even (respectively, odd) if each part $a_i$ is even (respectively, odd).
The empty composition $()$ is the only composition that is both even and odd.

We describe $\Pi_{+}$ in explicit terms.

\begin{prop}\label{P:E} Assume that $\ch(\F)\neq 2$.
The set $\{M_\alpha\}_{\alpha \text{ even}}$ is a linear basis for the even subalgebra $\Pi_{+}$
of
$(\QSym,\zetaQ)$.
\end{prop}
\begin{proof} Recall that $(\zetaQ)_n=H_n\in\NSym$, where $H_\alpha=M_\alpha^*$. We have 
\[(\BzetaQ)_n-(\zetaQ)_n=\begin{cases}
2H_n & \text{ if  $n$ is even,}\\
0 & \text{ if $n$ is odd.}
\end{cases}\]
Let $I$ be the ideal of $\NSym$ generated by $\{H_n\}_{\text{$n$ odd}}$. Then,
by~\eqref{E:prodH}, $\{H_\alpha\}_{\alpha \text{ not even}}$ is a linear basis of $I$, and by
Theorem~\ref{T:S-I}, 
\[\Pi_{+}=I^\perp=\{h\in\QSym\ \mid\ f(h)=0\ \ \forall f\in I\}\,.\]
Thus, $\{M_\alpha\}_{\alpha \text{ even}}$ is a linear basis for $\Pi_{+}$. 
\end{proof}

We now turn to $\Pi_{-}$. An explicit description in terms of linear relations (the generalized
Dehn-Sommerville relations) has already been given in Example~\ref{Ex:BB}. We  proceed to
describe a linear basis. We make use of the morphisms $\Theta$ and $\Psi_2:\QSym\to\QSym$
introduced in Examples~\ref{Ex:powerzeta}--\ref{Ex:theta}.

\begin{prop}\label{P:Theta-Psi2} Assume that $\ch(\F)\neq 2$. The morphisms
$\Theta:\QSym\to\QSym$ and
$\Psi_2:\QSym\to\QSym$ satisfy
\begin{itemize}
\item[(a)] $\Psi_2|_{\Pi_{-}}=\Theta|_{\Pi_{-}}$,
\item[(b)] $\Theta(\Pi_{-})=\Pi_{-}$,
\item[(c)] $\Theta(\QSym)=\Pi_{-}$.
\end{itemize}
\end{prop}
\begin{proof} By construction, $\Psi_2$ and $\Theta$ are the morphisms corresponding to the characters
$\zetaQ^2$ and $\nu=\BzetaQ^{-1}\zetaQ$, respectively. Now, by Definition~\ref{D:evenodd},
\[\BzetaQ^{-1}|_{\Pi_{-}}=\zetaQ|_{\Pi_{-}}\,.\]
Therefore, $\nu|_{\Pi_{-}}=\zetaQ^2|_{\Pi_{-}}$ and by uniqueness $\Theta|_{\Pi_{-}}=\Psi_2|_{\Pi_{-}}$.

Since $\nu$ is an odd character on $\QSym$, we have, by Corollary~\ref{C:E-Pi-univ},
\[\Theta(\QSym)\subseteq\Pi_{-}\,.\]
In particular, $\Theta(\Pi_{-})\subseteq\Pi_{-}$. On the other hand, formula~\eqref{E:powerpsi} shows
that $\Psi_2:\QSym\to\QSym$ is bijective (since $\ch(k)\neq 2$),  so by the above $\Theta|_{\Pi_{-}}$
is injective. Hence $\Theta(\Pi_{-})=\Pi_{-}$ and then also
$\Theta(\QSym)=\Pi_{-}$.
\end{proof}

Given an odd composition  $\beta$, set
\begin{equation}\label{E:eta}
\eta_{\beta}\ :=\ \sum_{\alpha\leq\beta} 2^{k(\alpha)}M_\alpha\,.
\end{equation}
The sum is over all compositions (not necessarily odd) smaller than $\beta$. Up to a sign, these
are the elements introduced by Hsiao in~\cite[Section 2]{Hsiao}.

\begin{prop}\label{P:Pibasis} Assume that $\ch(\F)\neq 2$. The odd subalgebra $\Pi_{-}$
of $(\QSym,\zetaQ)$ is Stembridge's Hopf algebra.
The set $\{\eta_\beta\}_{\beta \text{  odd}}$ is a linear basis for $\Pi_{-}$. 
\end{prop}
\begin{proof} In Example~\ref{Ex:theta} we identified $\Theta$ with the
morphism defined by Stembridge in~\cite{Ste97}. 
Stembridge's Hopf algebra is defined  as the image of $\Theta$.
By Proposition~\ref{P:Theta-Psi2}, this image is $\Pi_{-}$.

By triangularity, the set $\{\eta_\beta\}_{\beta \text{  odd}}$ is linearly
independent. On the other hand, from the explicit form of $\Theta$~\eqref{E:theta}, we see that
this set spans $\im(\Theta)=\Pi_{-}$. 
\end{proof}

\begin{rem}\label{R:Pi}
 The dimension of the $n$-th homogeneous component of $\Pi_{+}$, i.e., the number
of even
 compositions of $n$, is $2^{\frac{n}{2}-1}$ if $n>0$ is even and $0$ if $n$ is odd.
 The dimension of the $n$-th homogeneous component of $\Pi_{-}$, i.e., the number of odd compositions
of $n$, is the Fibonacci number $f_{n}$, where $f_1=f_2=1$ and
$f_n=f_{n-1}+f_{n-2}$ for $n\geq 3$.
\end{rem}
\begin{rem} We provide a description for the dual of $\Pi_{-}$. According to
Proposition~\ref{P:even-odd-basic},
\[\Pi_{-}^*=\NSym/I\,,\]
where $I:=I(\BzetaQ,\zetaQ^{-1})=I(\chiQ,\epsilonQ)$ is the ideal of $\NSym$ generated by the
elements
\[(\BzetaQ)_n-(\zetaQ^{-1})_n=(-1)^nH_n-\sum_{\alpha\comp n}(-1)^{k(\alpha)}H_{\alpha}\,,\ 
n>0\,,\] or also by the elements
\[(\chiQ)_n=(\chiQ)_n-(\epsilonQ)_n=\sum_{i=0}^n (-1)^i H_{(i,n-i)}\,,\  n>0\,.\]
Billera and Liu have shown that, if $\ch(\F)\neq 2$, $I$ is in fact generated by the elements
$(\chiQ)_n$ for $n$ even~\cite[Propositions 3.2, 3.3]{BL}. {}From our point of view, this may be
 seen as follows. Consider the linear functional
\[\omegaQ:=\zetaQ\chiQ+\chiQ\BzetaQ=\zetaQ\BzetaQ\zetaQ+\BzetaQ\zetaQ\BzetaQ\,.\]
Then $\Bar{\omega}_\calQ=\BzetaQ\zetaQ\BzetaQ+\zetaQ\BzetaQ\zetaQ=\omegaQ$, so
$(\omegaQ)_{n}=0$ for $n$ odd. Hence, if $n$ is odd,
\[0=2(\chiQ)_n+\sum_{i=1}^n (\zetaQ)_i(\chiQ)_{n-i}+\sum_{i=0}^{n-1}(\chiQ)_i(\zetaQ)_{n-i}\,.\]
It follows by induction that the elements $(\chiQ)_n$ for odd $n$ belong to the ideal generated
by the elements $(\chiQ)_n$ for even $n$.
\end{rem}

\subsection*{Canonical projections onto the odd and even subalgebras}\label{S:canonicalproj}

Consider the canonical decomposition of Theorem~\ref{T:decomposition} applied to the
character $\zetaQ$:
\[\zetaQ=(\zetaQ)_+(\zetaQ)_-\]
with $(\zetaQ)_+$ even and $(\zetaQ)_+$ odd. Let
\[\Psi_+:(\QSym,(\zetaQ)_+)\to(\QSym,\zetaQ) \text{ \ and \
}\Psi_-:(\QSym,(\zetaQ)_-)\to(\QSym,\zetaQ)\] 
be the unique morphisms of combinatorial Hopf algebras (Theorem~\ref{T:univ-qsym}). By
Corollary~\ref{C:E-Pi-univ},
\begin{equation}\label{E:+-inclusion}
\im(\Psi_+)\subseteq\Pi_{+} \text{ \ and \
}\im(\Psi_-)\subseteq\Pi_{-}\,.
\end{equation}
 We can in fact show that $\Psi_+$ and $\Psi_-$ are projections
onto $\Pi_{+}$ and $\Pi_{-}$:

\begin{prop}\label{P:canonicalproj} With the above notation,
\begin{itemize}
\item[(a)] $\im(\Psi_+)=\Pi_{+}$  and $\im(\Psi_-)=\Pi_{-}$.
\item[(b)] $\Psi_+|_{\Pi_{+}}=\id|_{\Pi_{+}}$ and $\Psi_-|_{\Pi_{-}}=
\id|_{\Pi_{-}}$.
\item[(c)] $\Psi_+|_{\Pi_{-}}=\epsilon|_{\Pi_{-}}$ and $\Psi_-|_{\Pi_{+}}
=\epsilon|_{\Pi_{+}}$.
\item[(d)] $\Psi_+\Psi_-=\id$ (convolution product in $\End(\QSym)$).
\end{itemize}
\end{prop}
\begin{proof} By definition of even and odd subalgebras (Definition~\ref{D:evenodd}), 
$\zetaQ|_{\Pi_{+}}$ is even and $\zetaQ|_{\Pi_{-}}$ is odd. Hence, by uniqueness of the decomposition
of Theorem~\ref{T:decomposition}, 
\[\zetaQ|_{\Pi_{+}}=(\zetaQ)_+|_{\Pi_{+}}\,,\ \ \zetaQ|_{\Pi_{-}}=(\zetaQ)_-|_{\Pi_{-}}\,,
\text{ \ and \ }
(\zetaQ)_+|_{\Pi_{-}}=\epsilon|_{\Pi_{-}}\,,\ \ (\zetaQ)_-|_{\Pi_{+}}=\epsilon|_{\Pi_{+}}\,.\]
Assertions (b) and (c) follow from here, in view of uniqueness in Theorem~\ref{T:univ-qsym}.
Together with~\eqref{E:+-inclusion}, (b) implies (a). 

To obtain (d) it suffices to show that $\zetaQ\circ(\Psi_+\Psi_-)=\zetaQ$,  
again by Theorem~\ref{T:univ-qsym}. Now, since $\zetaQ$ is a morphism of 
algebras, precomposing with $\zetaQ$ preserves convolution products. Hence,
\[\zetaQ\circ(\Psi_+\Psi_-)=(\zetaQ\circ\Psi_+)(\zetaQ\circ\Psi_-)=
(\zetaQ)_+(\zetaQ)_-=\zetaQ\,,\]
by definition of $\Psi_{\pm}$ and $(\zetaQ)_{\pm}$.
The proof is complete.
\end{proof}

\begin{prop}\label{P:Theta-Psi-} The map $\Theta$ and the projection 
$\Psi_-$ are related as follows.
\begin{itemize}
\item[(a)] $\Theta\circ\Psi_-=\Theta=\Psi_-\circ\Theta$.
\item[(b)] $\im(\Psi_-)=\im(\Theta)=\Pi_{-}$.
\item[(c)] $\ker(\Psi_-)=\ker(\Theta)$.
\end{itemize}
\end{prop}
\begin{proof} Assertion (b) is part of Propositions~\ref{P:Theta-Psi2} 
and~\ref{P:canonicalproj}. The fact that $\Psi_-\circ\Theta=\Theta$ also
follows from those propositions, since 
$\Psi_-|_{\im(\Theta)}=\id|_{\im(\Theta)}$.

By our construction of $\Theta$ and  Remark~\ref{R:squares},
\[\zetaQ\circ\Theta=\nuQ=\bigl((\zetaQ)_-\bigr)^2\,.\]
 Hence,
\[\zetaQ\circ\Theta\circ\Psi_-=\bigl((\zetaQ)_-\bigr)^2\circ\Psi_-=
\bigl((\zetaQ)_-\circ\Psi_-\bigr)^2=(\zetaQ\circ\Psi_-)^2=\bigl((\zetaQ)_-\bigr)^2\,,\]
 since $\im(\Psi_-)=\Pi_{-}$ and $(\zetaQ)_-|_{\Pi_{-}}=\zetaQ|_{\Pi_{-}}$.

Thus, $\zetaQ\circ\Theta=\zetaQ\circ\Theta\circ\Psi_-$ and by uniqueness in
Theorem~\ref{T:univ-qsym}, $\Theta\circ\Psi_-=\Theta$. This completes the proof
of (a).

{}From (a) it follows that $\ker(\Psi_-)\subseteq\ker(\Theta)$, which with (b) 
implies (c).
\end{proof}

We conclude with an important result that states that $\QSym$ is the
{\em product} of its even and odd subalgebras in the category of graded connected Hopf
algebras.

\begin{theo}\label{T:product} 
Let $\H$ be a graded connected Hopf algebra and 
\[\alpha_+:\H\to\Pi_+ \text{ \ and \ }\alpha_-:\H\to\Pi_-\]
two morphisms of graded Hopf algebras. Then there is a unique
morphism of graded Hopf algebras $\alpha:\H\to\QSym$ such that
the following diagram commutes:
\begin{equation}\label{E:product}
\raisebox{25pt}{\xymatrix{
& {\H}\ar[ld]_{\alpha_-}\ar@{-->}[d]^{\alpha}\ar[rd]^{\alpha_+} &\\
{\Pi_-} & {\QSym}\ar[l]^{\Psi_-}\ar[r]_{\Psi_+} & {\Pi_+}
}}
\end{equation}
\end{theo}
\begin{proof} Define two characters on $\H$ by
\[\zeta_+:=\zetaQ\circ\alpha_+ \text{ \ and \ }\zeta_-:=\zetaQ\circ\alpha_-\,.\]
In other words, $\zeta_{\pm}=\alpha^*_{\pm}(\zetaQ|_{\Pi_{\pm}})$.
Since $\zetaQ|_{\Pi_{+}}$ is even and $\zetaQ|_{\Pi_{-}}$ is odd, $\zeta_+$ is
an even character of $\H$ and $\zeta_-$ is an odd character of $\H$. Define a
third character on $\H$ by
\[\zeta:=\zeta_+\zeta_-\,.\]
Let $\alpha:\H\to\QSym$ be the unique morphism of graded Hopf algebras such that
$\zetaQ\circ\alpha=\zeta$ (Theorem~\ref{T:univ-qsym}). Then, we have
\[\zeta=\alpha^*(\zetaQ)=\alpha^*\bigl((\zetaQ)_+(\zetaQ)_-\bigr)=
\alpha^*\bigl((\zetaQ)_+\bigr)\alpha^*\bigl((\zetaQ)_-\bigr)\,.\]
This is another decomposition of $\zeta$ as a product of an even and an odd
characters. By Theorem~\ref{T:decomposition} we must have
\[\zeta_{\pm}=\alpha^*\bigl((\zetaQ)_{\pm}\bigr)=(\zetaQ)_{\pm}\circ\alpha\,.\]
On the other hand, by definition of $\Psi_+$ and $\Psi_-$,
$(\zetaQ)_{\pm}=\zetaQ\circ\Psi_{\pm}$. Therefore,
\[\zetaQ\circ\Psi_{\pm}\circ\alpha=(\zetaQ)_{\pm}\circ\alpha=\zeta_{\pm}=
\zetaQ\circ\alpha_{\pm}\,.\]
It follows from Theorem~\ref{T:univ-qsym} that 
$\Psi_{\pm}\circ\alpha=\alpha_{\pm}$, i.e., diagram~\eqref{E:product} commutes.

Suppose $\tilde{\alpha}:\H\to\QSym$ is another morphism of graded Hopf
algebras such that $\Psi_{\pm}\circ\tilde{\alpha}=\alpha_{\pm}$.
In the derivation below, we make use of the fact that $\Psi_+\Psi_-=\id$
(Proposition~\ref{P:canonicalproj}):
\begin{align*}
\zetaQ\circ\tilde{\alpha}&=\zetaQ\circ(\Psi_+\Psi_-)\circ\tilde{\alpha}=
\zetaQ\circ\Bigl((\Psi_+\circ\tilde{\alpha})(\Psi_-\circ\tilde{\alpha})\Bigr)\\
&=\zetaQ\circ(\alpha_+\alpha_-)=(\zetaQ\circ\alpha_+)(\zetaQ\circ\alpha_-)=
\zeta_+\zeta_-=\zeta=\zetaQ\circ\alpha\,.
\end{align*}
{}From Theorem~\ref{T:univ-qsym} we deduce $\tilde{\alpha}=\alpha$.
\end{proof}

\begin{rem} The {\em Hilbert series} (see Section~\ref{S:other}) 
of $\QSym$, $\Pi_+$, and $\Pi_-$ are
\[\QSym(x)=\frac{1-x}{1-2x}\,,\ \ \ \Pi_+(x)=\frac{1-x^2}{1-2x^2}\,,
 \text{ \ \ and \ \ }\Pi_-(x)=\frac{1-x^2}{1-x-x^2}\,.\]
The latter two expressions follow from Propositions~\ref{P:E} 
and~\ref{P:Pibasis}.
Note that
\[\frac{1}{\QSym(x)}=\frac{1}{\Pi_+(x)}+\frac{1}{\Pi_-(x)}-1\,.\]
This may also be seen as a consequence of Theorem~\ref{T:product}. In fact,
it follows from this theorem that $\NSym=(\QSym)^*$ is the {\em free
product} of the algebras $(\Pi_+)^*$ and $(\Pi_-)^*$, and it is a general
fact that the Hilbert series of a free product $A=A_1\ast
A_2$ of two graded connected algebras is related to the Hilbert series of 
the factors by means of
\begin{equation}\label{E:freeproduct}
\frac{1}{A(x)}=\frac{1}{A_1(x)}+\frac{1}{A_2(x)}-1\,.
\end{equation}
\end{rem}
 %%%%%%%%%%%%%%%%%%%%%%%%%%%%%%%%%%%%%%%%%%%%%%%%%%%%%
 \section{Even and odd subalgebras of $\Sym$}\label{S:even-odd-Sym}

We describe the even and odd subalgebras of the combinatorial Hopf algebra of
symmetric functions $(\Sym,\zetaS)$. We find that the odd subalgebra is the Hopf algebra of
Schur $Q$-functions~\cite[III.8]{Mac}.

Let $\lambda=(l_1,\ldots,l_k)$ be a partition of $n$. As for compositions, we write
$\abs{\lambda}=n$ and
$k(\lambda)=k$. Recall the power sum symmetric functions
$p_\lambda$, defined by
\begin{equation}\label{E:prodp}
 p_n\ :=\ m_n\text{ \ and \ }
 p_\lambda\ :=\ p_{l_1}\cdots p_{l_k}\,.
\end{equation}

As for compositions, we say that a partition is even (odd) if all its parts are even (odd).
For each odd partition $\lambda$, define an element
\begin{equation}\label{E:eta-symm}
\eta_\lambda\ :=\ \sum_{s(\alpha)=\lambda}\eta_\alpha\,,
\end{equation}
where $\eta_\alpha$ is as in~\eqref{E:eta} and $s(\alpha)$ as in~\eqref{E:monomial}.

\begin{prop}\label{P:odd-Sym} Assume $\F=\QQ$. The sets $\{p_\lambda\}_{\lambda \text{ odd}}$
and  $\{\eta_\lambda\}_{\lambda \text{ odd}}$ are linear bases for the odd subalgebra of
$(\Sym,\zetaS)$.
\end{prop}
\begin{proof} Consider first the elements $\eta_\lambda$. By
item (f) in Proposition~\ref{P:even-odd-basic} applied to $\Sym\inc\QSym$, 
\[S_-(\Sym,\zetaS)=S_-(\QSym,\zetaQ)\cap\Sym\,.\]
The set $\{\eta_\alpha\}_{\alpha \text{ odd}}$ is a basis for
$S_-(\QSym,\zetaQ)$, by Proposition~\ref{P:Pibasis}. It follows readily
from~\eqref{E:monomial}, ~\eqref{E:eta}, and~\eqref{E:eta-symm} that the set
$\{\eta_\lambda\}_{\lambda \text{ odd}}$ is a basis for $S_-(\Sym,\zetaS)$.

Note that $\eta_{2n+1}=m_{2n+1}=p_{2n+1}$. Hence, by~\eqref{E:prodp}, and since
$S_-(\Sym,\zetaS)$ is a subalgebra, $p_\lambda\in S_-(\Sym,\zetaS)$ for each odd partition
$\lambda$. Since the set  $\{p_\lambda\}_{\lambda \text{ odd}}$ is
 linearly independent and has the same number of elements in each homogeneous dimension as the
set $\{\eta_\lambda\}_{\lambda \text{ odd}}$,  it must also be a linear
basis for $S_-(\Sym,\zetaS)$.
\end{proof}

\begin{rem} It is easy to see that 
\[\zetaS(p_\lambda)=1\,,\ \ \BzetaS(p_\lambda)=(-1)^{\abs{\lambda}}\,, \text{ \ and \ }
\zetaS^{-1}(p_\lambda)=(-1)^{k(\lambda)}\,.\]
One may use these to show directly that $\{p_\lambda\}_{\lambda \text{ odd}}$ is a basis 
for $S_-(\Sym,\zetaS)$.
\end{rem}

\begin{prop}\label{P:even-Sym}Assume $\F=\QQ$. The sets $\{p_\lambda\}_{\lambda \text{ even}}$
and
$\{m_\lambda\}_{\lambda \text{ even}}$ are linear bases
for the even subalgebra of $(\Sym,\zetaS)$.
\end{prop}
\begin{proof} By item (f) in Proposition~\ref{P:even-odd-basic} applied to $\Sym\inc\QSym$, 
\[S_+(\Sym,\zetaS)=S_+(\QSym,\zetaQ)\cap\Sym\,.\]
According to Proposition~\ref{P:E}, the set $\{M_\alpha\}_{\alpha \text{ even}}$ is a linear
basis for $S_+(\QSym,\zetaQ)$. It is then clear from~\eqref{E:monomial} that the set
$\{m_\lambda\}_{\lambda \text{ even}}$ is a basis for $S_+(\Sym,\zetaS)$. Recall that
$p_{2n}=m_{2n}$. As in the proof of Proposition~\ref{P:odd-Sym}, it follows that the set
$\{p_\lambda\}_{\lambda \text{ even}}$ is a linear basis for $S_+(\Sym,\zetaS)$.
\end{proof}

\begin{rem}  Proposition~\ref{P:odd-Sym} shows that the odd subalgebra of $(\Sym,\zetaS)$ is
the Hopf algebra of Schur $Q$-functions, as defined in~\cite[III.8]{Mac}.

 As seen in the above proof, $S_-(\Sym,\zetaS)=S_-(\QSym,\zetaQ)\cap\Sym$.
This says that (if $\F=\QQ$) the intersection of Stembridge's Hopf algebra $\Pi_{-}$ with $\Sym$ is
the Hopf algebra of Schur $Q$-functions, a known result from~\cite[Theorem 3.8]{Ste97}.

 Recall that the graded Hopf algebra $\Sym$ is self-dual~\eqref{E:mhduality}. Thus, the
components of a character on $\Sym$ are elements of
$\Sym$ itself, and Proposition~\ref{P:duality} relates even and odd characters on $\Sym$ to
the even and odd subalgebras of $(\Sym,\zetaS)$.
Consider the canonical odd character $\nuS$ and the Euler character $\chiS$. Since
$\Sym$ is cocommutative,  $\chiS$ is even (Example~\ref{Ex:nu-chi}).
Thus, Proposition~\ref{P:duality} guarantees that the components of $\nuS$ belong to
$S_-(\Sym,\zetaS)$ and those of $\chiS$ to $S_+(\Sym,\zetaS)$.  Let $e_n$ be the 
elementary symmetric functions. The components of the canonical characters of $\Sym$ are:
\begin{align*}
(\zetaS)_n & = h_n=\sum_{\lambda\parti n}m_\lambda\,,\\
(\zetaS^{-1})_n & =(-1)^ne_n\,,\\
(\nuS)_n & =\sum_{i=0}^n e_i h_{n-i}\,,\\
(\chiS)_n & = \sum_{i=0}^n (-1)^i h_{(i,n-i)}=\begin{cases} 
\sumsub{\lambda\parti n\\\lambda\text{ even}}m_\lambda & \text{ if $n$ is even,}\\
0 & \text{ if $n$ is odd.}
\end{cases}
\end{align*}
In particular, we see that the components of $\nuS$ are precisely the Schur
 $Q$-functions $q_n$, as defined in~\cite[Chapter III.8, formula (8.1)]{Mac}.  
\end{rem}

\begin{rem} The decomposition of the Hopf algebra of quasi-symmetric functions
as the categorical product of its even and odd Hopf subalgebras (Theorem~\ref{T:product})
admits an analog for symmetric functions:
$\Sym$ is the product of its even and odd subalgebras
in the category of cocommutative Hopf algebras. 
In this category, the product is simply the tensor product of Hopf algebras.
Thus,
\begin{equation}\label{E:sym-product}
\Sym\cong S_+(\Sym,\zetaS)\otimes S_-(\Sym,\zetaS)\,.
\end{equation}
The proof of this fact is completely analogous to that of 
Theorem~\ref{T:product}. It also follows at once from the explicit
descriptions for $S_{\pm}(\Sym,\zetaS)$ in Propositions~\ref{P:odd-Sym} and~\ref{P:even-Sym}.

In terms of Hilbert series, the decomposition~\eqref{E:sym-product} 
boils down to
\[\prod_{n\geq 1}\frac{1}{1-x^{n}}\ =\ 
\prod_{n\geq 1}\frac{1}{1-x^{2n}}\prod_{n\geq 1}\frac{1}{1-x^{2n-1}}\,.\]
Since $\Sym$ is commutative, \eqref{E:sym-product} is also a {\em coproduct} decomposition
in the category of commutative Hopf algebras. In fact, this decomposition is
self-dual.
\end{rem}

 %%%%%%%%%%%%%%%%%%%%%%%%%%%%%%%%%%%%%%%%%%%%%%%%%%%%%%%%%%%%%%%%%%
 \section{The odd subalgebras of other combinatorial Hopf algebras}\label{S:other}

Symmetric, quasi-symmetric, and non-commutative symmetric functions
(Section~\ref{S:Q-N}) are related by the following commutative diagram of graded Hopf algebras:
\begin{equation}\label{E:CD}
\raisebox{-25pt}{ \begin{picture}(220,65)(0,-25) 
 \put(0,0){$\Sym$} \put(120,0){$\YSym$} \put(190,0){$\SSym$} 
\put(55,-25){$\NSym$} \put(55,25){$\QSym$} 
\put(23,  9){\epsffile{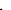}} 
\put(23,  9){\vector(2,1){30}} 
 \put(53,-18){\vector(-2,1){30}} 
 \put(53,-18){\vector(-2,1){25}} 
\put( 88,-18){\epsffile{dhook.eps}} 
\put( 88,-18){\vector(2,1){30}} 
\put(118,  9){\vector(-2,1){30}} 
\put(118,  9){\vector(-2,1){25}} 
 \put(158,5){\epsffile{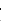}} 
\put(158,5){\vector(1,0){27}} 
 \put(185,0){\vector(-1,0){30}}
\put(185,0){\vector(-1,0){25}}
 \end{picture} }
\end{equation}

The Hopf algebra $\SSym$ has linear bases indexed by all permutations of $[n]$, $n\geq 0$, and
was introduced by Malvenuto and Reutenauer~\cite{Malv,MR}. The Hopf algebra
$\YSym$ has linear bases indexed by planar binary trees and was introduced by Loday and
Ronco~\cite{LR,LR02}. See also~\cite{AS,ASb} for a detailed study of these Hopf algebras.
Both $\SSym$ and $\YSym$ are self-dual Hopf algebras (in the graded sense). In
diagram~\eqref{E:CD}, duality acts as reflection across the middle horizontal line (not just on
objects but also on morphisms).

In Sections~\ref{S:even-odd-QSym} and~\ref{S:even-odd-Sym} we calculated the even and odd
subalgebras of $\QSym$ and $Sym$. In this section we calculate the odd subalgebras of $\NSym$,
$\SSym$ and $\YSym$. The even subalgebras may be calculated by the same methods. We view these as
combinatorial Hopf algebras via the pull-back of the canonical character $\zetaQ$ of $\QSym$;
therefore, \eqref{E:CD} is a commutative diagram of morphisms of combinatorial
Hopf algebras.

In this section, we use $\OQ$ to denote the odd subalgebra of $(\QSym,\zetaQ)$
(instead of $S_-(\QSym,\zetaQ)$), and $\IQ$ to denote the odd ideal $I_-(\QSym,\zetaQ)$.
Thus, $\OQ$ is a Hopf subalgebra of $\QSym$, $\IQ$ is a Hopf ideal of $\NSym$, and by
Proposition~\ref{P:even-odd-basic}\,(a),
\[(\OQ)^*=\NSym/\IQ\,.\]
Similarly, we use $\OS$ and $\IS$, $\ON$ and $\IN$, $\OSS$ and $\ISS$, and $\OY$ and
$\IY$ for the odd subalgebras and ideals associated to the combinatorial Hopf algebras $\Sym$,
$\NSym$, $\YSym$, and $\SSym$.

If $V=\bigoplus_{n\geq 0}V_n$ is a graded vector space and each $V_n$ is finite-dimensional,
we let
\[V(x):=\sum_{n\geq 0}\dim V_n\,x^n\]
denote the corresponding Hilbert series.

{}From Proposition~\ref{P:Pibasis} and Remark~\ref{R:Pi} we know that 
\begin{equation}\label{E:odd-qsym}
\OQ(x)=\frac{1-x^2}{1-x-x^2}\,.
\end{equation}

{}From Proposition~\ref{P:odd-Sym} it follows that 
\begin{equation}\label{E:odd-sym}
\OS(x)=\prod_{n\geq 0}\frac{1}{1-x^{2n+1}}\,.
\end{equation}

Below we provide formulas for the Hilbert series of $\ON$, $\OSS$, and $\OY$.

\subsection*{The odd subalgebra of $\NSym$} We describe the graded dual of
$\ON$. This is a quotient Hopf algebra of $(\NSym)^*=\QSym$.

Our derivation is based on the following simple fact. Let $V$ be a vector space and
$S(V)$ the free commutative algebra on $V$ (the symmetric algebra). Let $W$ be another vector
space and consider the inclusion
$S(V)\inc S(V\oplus W)$ induced by the canonical map $V\inc V\oplus W$.
Let $I$ be an ideal of $S(V)$ and $J$ the ideal of $S(V\oplus W)$ generated by $I$.
Then there are canonical isomorphisms of algebras
\begin{equation}\label{E:ideal-comm}
\frac{S(V\oplus W)}{J}=\frac{S(V)\otimes S(W)}{I\otimes S(W)}=\frac{S(V)}{I}\otimes
S(W)\,.
\end{equation}

The algebra of symmetric functions is freely generated, as a commutative algebra, by the
elements $m_n$, $n\geq 1$~\cite[Chapter I.2]{Mac}. It is also known that the algebra of
quasi-symmetric functions is freely generated, as a commutative algebra over $\QQ$, by the
elements $M_\alpha$ indexed by {\em Lyndon} compositions $\alpha$~\cite[Theorem 2.6, Example 6.1]{Ho00}. 
The latter set of generators contains the former, as all compositions with one part are Lyndon
and $M_{(n)}=m_n$. (A similar result is~\cite[Corollary 2.2]{MR}.)

Let $V$ be the vector subspace of $\Sym$ spanned by $\{m_n\}_{n\geq 1}$ and
$W$ the subspace of $\QSym$ spanned by the elements $M_\alpha$ where $\alpha$ runs
over the set of Lyndon compositions with $2$ or more parts. Thus, as algebras,
\begin{equation}\label{E:QSym-over-Sym}
\Sym=S(V),\ \ \QSym=S(V\oplus W)=S(V)\otimes S(W)\,,
\end{equation}
and the inclusion $\Sym\inc\QSym$ is induced by the canonical map $V\inc V\oplus W$.

Since the map $\NSym\to\Sym$ is
a morphism of combinatorial Hopf algebras, $\IN$ is
the ideal of $(\NSym)^*=\QSym$ generated by $\IS$ (by item (e) in
Proposition~\ref{P:even-odd-basic}). Therefore, applying~\eqref{E:ideal-comm} with 
$I=\IS$ and $J=\IN$ and making use of
 item (a) in Proposition~\ref{P:even-odd-basic} we obtain
\begin{equation}\label{E:odd-nsym-dual}
(\ON)^*=\frac{\QSym}{\IN}=\frac{\Sym}{\IS}\otimes S(W)=(\OS)^*\otimes S(W)
\end{equation}
(we used $(\NSym)^*=\QSym$ and $(\Sym)^*=\Sym$). 

For $n\geq 1$, the number of Lyndon compositions of $n$ is~\cite[Sequence A059966]{Sl}
\[\ell_n:=\frac{1}{n}\sum_{d\mid n}\mu\Bigl(\frac{n}{d}\Bigr)(2^d-1)\,,\]
where $\mu$ is the classical (number-theoretic) M\"obius function. Therefore,
$S(W)(x)=\prod_{n\geq 1}\frac{1}{(1-x^n)^{\ell_n-1}}$, and from~\eqref{E:odd-sym}
and~\eqref{E:odd-nsym-dual} we deduce
\[\ON(x)=\prod_{n\geq 1}\frac{1}{(1-x^n)^{a_n}}\,,
\text{ \ where \ } a_n=\begin{cases}\ell_n & \text{ if $n$ is odd,}\\ \ell_n-1 & \text{ if $n$
is even.}\end{cases}\]
Thus,
\[\ON(x)=1+x+x^2+3x^3+5x^4+11x^5+22x^6+44x^7+\cdots\,.\]
Note that~\eqref{E:QSym-over-Sym} and~\eqref{E:odd-nsym-dual} also imply that
\[\frac{\ON(x)}{\OS(x)}=\frac{\NSym(x)}{\Sym(x)}\,.\]

\subsection*{The odd subalgebras of $\SSym$ and $\YSym$} We describe the graded dual of
$\OSS$ by a derivation similar to the one in the previous section. 

Let $V$ be a vector space and $T(V)$ the free algebra on $V$ (the tensor algebra). Let $W$ be
another vector space and consider the inclusion
$T(V)\inc T(V\oplus W)$ induced by the canonical map $V\inc V\oplus W$.
Let $I$ be an ideal of $T(V)$ and $J$ the ideal of $T(V\oplus W)$ generated by $I$.
There is a canonical isomorphism of algebras
\begin{equation}\label{E:ideal-noncomm}
\frac{T(V\oplus W)}{J}=\frac{T(V)}{I}\ast T(W)\,,
\end{equation}
where $\ast$ denotes the free product of algebras.

The algebra of non-commutative symmetric functions is freely generated 
by the elements $H_n$, $n\geq 1$. According to~\cite[Theorem 6.1, Remark 6.2]{AS}, 
the algebra $(\SSym)^*$ is freely generated by certain elements indexed by permutations with no
{\em global descents}. Moreover, it follows from~\cite[Theorem 7.3]{AS} that the inclusion
$\NSym\inc (\SSym)^*$ sends $H_n$ to the basis element indexed by the identity permutation 
of $n$ (which has no global descents). Thus, there are spaces $V$ and $W$ such that
\begin{equation}\label{E:SSym-over-NSym}
\NSym=T(V)\,,\ \  (\SSym)^*=T(V\oplus W)=T(V)\ast T(W)\,, 
\end{equation}
and the inclusion $\NSym\inc (\SSym)^*$ is induced by $V\inc V\oplus W$. The space $W$ has a
basis indexed by permutations that have no global descents and are different from the identity.

Since the map $\SSym\to\QSym$ is
a morphism of combinatorial Hopf algebras, $\ISS$ is
the ideal of $(\SSym)^*$ generated by $\IQ$ (by item (e) in
Proposition~\ref{P:even-odd-basic}). Therefore, applying~\eqref{E:ideal-noncomm} with 
$I=\IQ$ and $J=\ISS$  we obtain
\begin{equation}\label{E:odd-ssym-dual}
(\OSS)^*=\frac{(\SSym)^*}{\ISS}=\frac{\NSym}{\IQ}\ast T(W)=(\OQ)^*\ast T(W)\,.
\end{equation}

Recall that the Hilbert series of a free product is related to 
the Hilbert series of the factors by means of~\eqref{E:freeproduct}.
Applying this to~\eqref{E:SSym-over-NSym} and~\eqref{E:odd-ssym-dual} we deduce 
\[\frac{1}{\OSS(x)}-\frac{1}{\OQ(x)}=\frac{1}{\SSym(x)}-\frac{1}{\QSym(x)}\,.\]

This reduces to
\[\frac{1}{\OSS(x)} =\frac{1}{\sum_{n\geq 0}n!\,x^n} +\frac{x^2}{1-x^2}\,,\]
from which
\[\OSS(x)=1+x+x^2+4x^3+19x^4+105x^5+660x^6+\cdots\,.\]

\medskip 
The Hopf algebra $(\YSym)^*$ is also free on a set of generators containing the generators of
$\NSym$, according to the results of~\cite{ASb}. Therefore analogous considerations apply and
we obtain
\[\frac{1}{\OY(x)}-\frac{1}{\OQ(x)}=\frac{1}{\YSym(x)}-\frac{1}{\QSym(x)}\,.\]

 %%%%%%%%%%%%%%%%%%%%%%%%%%%%%%%%%%%%%%%%%%%%%%%%%%%%%%%%%%%%%%

 \end{document}